\documentclass[11pt]{article}
\usepackage[dvipsnames]{xcolor}
\usepackage{pst-plot,pstricks,pst-node}
\usepackage[latin1]{inputenc}
\usepackage{xspace}
\usepackage{epsfig,graphicx}
\usepackage{amsmath,amssymb}
\usepackage{theorem}
\usepackage{setspace}
\usepackage{multirow}
\usepackage{authblk}
\usepackage{fancyhdr} 
\usepackage[numbers,sort]{natbib} 
\onehalfspace
\usepackage{mathtools}
\usepackage{tabls}

\usepackage[algo2e,lined,algoruled,titlenumbered]{algorithm2e}
\allowdisplaybreaks[1]

\graphicspath{pics/}
\DeclareGraphicsExtensions{.eps}

\fancypagestyle{memo}{
	\pagestyle{fancy}
	\fancyhf{}%
	\fancyhead[RO]{\thepage}%
	\fancyhead[LO]{\rightmark}%
	
}%

\theoremstyle{plain}
\newtheorem{theorem}{Theorem}
\newtheorem{prop}{Proposition}
\newtheorem{lemma}{Lemma}

\newtheorem{example}{Example}
\newtheorem{observation}{Observation}

\theoremstyle{plain} {
	\theorembodyfont{\rmfamily}%

}

\newcommand{\ProofNoNL}{{\bf \noindent Proof.}\xspace}

\newcommand{\EndProofNoNL}{\hfill $\Box$ \par \bigskip}

\newcommand{\ExampleCont}[1]{\vspace{2ex} \noindent {\bf Example #1 (cont.)}\itshape}
\newcommand{\EndExample}{\normalfont \vspace{2ex}}
\newcommand{\wloge}{without loss of generality\xspace}

\newcommand{\Real}{\mathbb{R}\xspace}

\newcommand{\ex}{\ensuremath{e_x}\xspace}
\newcommand{\ey}{\ensuremath{e_y}\xspace}
\newcommand{\tx}{\ensuremath{t_x}\xspace}
\newcommand{\name}{\mbox{(MMR-EMCLP)}\xspace}

\usepackage[normalem]{ulem}

\begin{document}

\title{Minmax Regret Maximal Covering Location Problems with Edge Demands}

\author[(a)]{Marta Baldomero-Naranjo\footnote{Corresponding author: marta.baldomero@uca.es}}
\author[(b)]{J\"org Kalcsics}
\author[(a)]{Antonio M. Rodr\'iguez-Ch\'ia}

\affil[(a)]{\small{Departamento de Estad\'istica e Investigaci\'on Operativa, Universidad de C\'adiz, C\'adiz, Spain, marta.baldomero@uca.es, antonio.rodriguezchia@uca.es}}
\affil[(b)]{\small{School of Mathematics, The University of Edinburgh, Edinburgh, United Kingdom, joerg.kalcsics@ed.ac.uk}}

\maketitle

\begin{abstract}
This paper addresses a version of the single-facility Maximal Covering Location Problem on a network where the demand is: i) distributed along the edges and ii) uncertain with only a known interval estimation. To deal with this problem, we propose a minmax regret model where the service facility can be located anywhere along the network. This problem is called Minmax Regret Maximal Covering Location Problem with demand distributed along the edges \name. Furthermore, we present two polynomial algorithms for finding the location that minimises the maximal regret assuming that the demand realisation is an unknown constant or linear function on each edge. We also include two illustrative examples as well as a computational study for the unknown constant demand case to illustrate the potential and limits of the proposed methodology.

\end{abstract}
\bigskip


\section{Introduction}%
\label{sec:Introduction}%
In this work, we focus our research on the Maximum Covering Location Problem (MCLP), first introduced by \citet{ChuRev74}. For a given set of demand points (clients or customers), the goal of this problem is to locate one or more service facilities in such a way that the demand of covered clients is maximized. A client is considered to be covered if the distance to an open facility is smaller than or equal to a fixed coverage radius, $R$. We can find many different variants of this problem in the literature, some of them in continuous space and others on networks or in discrete spaces (see \citet{Pla02,GarMar15,BerKra02} for a detailed survey of each case). In the past, almost all of the models discussed in the literature assume that the demand is represented by a finite set of points. However, recently, an increasing body of literature in facility location addresses problems where the demand is continuously distributed over a region (polygon, straight lines, edges, etc.), see \citet{MurTon07,MurTonKim10}.

In the  case of network problems (the subject of study of this paper), there are several applications where the assumption that the demand only occurs at the nodes of the network is not realistic; e.g.\ for the location of emergency facilities (see \citet{LiZhaZhuWya11}), for the planning of garbage collection, meter reading, or mail delivery, or situations where the coverage areas are extremely distance-dependent such as the location of Automated External Defibrillators (AED), bus stops, Automated Teller Machines (ATM), etc. Similarly, assuming the situation that the service facilities can only be located at the nodes of the network, the so-called node-restricted version, also leads to unrealistic models for some situations, see \citet{RevTorFal76,ChuMea79}. For example, in the case of locating bus stops, this assumption  would mean that they can only  be placed at the corner of streets. Thus, assuming that the demand is concentrated at nodes and/or service facilities can only be located at nodes may yield models that are far from the real situation under study, providing solutions that are of little use, see \citet{CurSchi90} for further details.  

In the following, we review the main models in the literature related to our problem. The MCLP where the demand is assumed to exist only at nodes is solved in \citet{ChuRev74}. They proposed an integer programming formulation to solve the multi-facility node-restricted problem. \citet{GarMar15} and  \citet{MarMartRodSal18} provided an overview of models for the discrete covering location problem. Recently, \citet{CorFurLju19} proposed an effective decomposition approach for the MCLP based on a branch-and-Benders-cut formulation to solve realistic cases where the number of customers is much larger than the number of potential facility locations. To solve the continuous version of the problem, i.e., facilities can be located not only at the nodes but also along the edges, \citet{ChuMea79} proposed a Finite Dominating Set (FDS). An FDS for a location problem is a finite set of points that is guaranteed to contain an optimal solution for the problem. If such a set can be found, then the continuous optimization problem can equivalently be formulated as a discrete optimization problem and solved applying techniques of discrete optimization. This technique is widely applied in the literature to solve these kinds of problems, see for example  \citet{NicPue99, KalNicPueTam02, KalNicPue03,BerKalKraNic09}.

The MCLP where the demand is distributed along the edges (edge demand) and the facilities can be located anywhere on the graph (continuous location space) was first proposed by \citet{BerKalKra16}. They found a FDS for the case of known demand functions and presented an exact approach to solve the single facility problem for constant and piecewise linear demand functions (unlike our paper, the actual demand functions are assumed to be known a priori). Extensions of this model are also analysed, as for instance the obnoxious version, the multi-facility case, and the conditional problem (when all but one facility have already been located). Moreover, they proposed a heuristic for the multi-facility case. Under the same framework (edge demand and continuous location space), a minimum set-covering problem (locating the minimum number of facilities to cover all the demand points) and a stochastic version of the MCLP are studied by \citet{FrMaHa20}  and \citet{BlaCarBog16}, respectively. In the later case, the expected demand covered is maximized and the problem is formulated as a Mixed Integer Nonlinear Program (MINLP) that is solved using a branch-and-bound algorithm comprising a combinatorial part (in which edges of the network are chosen to contain facilities) and a continuous global optimization part (once the edges are chosen, the optimal locations are found on those edges).

In this paper, we consider the MCLP under this framework on a network, i.e., edge demand and  continuous location space. However, unlike the aforementioned papers, the actual intensity of the demand along each edge is unknown. Instead, the only available knowledge is an upper and lower bound of this intensity at each point of the network.
We are not aware of any other paper dealing simultaneously with edge demand, a continuous location space, and uncertainty in the intensity of the demand for the MCLP on networks (or, in fact, for any location problem on networks). 
From a practical point of view, assuming that demand along an edge is known corresponds to an ideal but usually unrealistic scenario. Since demand is uncertain in nature and varies from one day to another, or even within a day, we will treat demands as being unknown. Nevertheless, we usually have a good  estimation of the minimal or maximal demand along the edge, so that we can at least assume demand to lie within a known range.

	Concerning real-life applications, the MCLP has been applied in several fields such as the location of emergency facilities, health care, natural disaster rescue, ecology, signal-transmission facilities, bike-sharing, police resources, see \citet{AR97,ARG02,Chung86,DD04,MuLiMuk20,WLN06,YanXiaZhanZhouYanXu20}. More concretely, our model can be used to locate AED, ATM, bus stops, automated parcel lockers, or bicycle parking racks, among others. In these applications the network represents a city, where the edges are the streets and the nodes are the intersections between them. Moreover, these services can be located anywhere along the network because they do not require sophisticated infrastructures or spaces. We consider that the demand is the potential amount of users for that service. Also, a client would be considered as a covered client if the distance between him/her and the service is lower than or equal to a fixed radius. In these applications, the exact demand of these services is usually unknown but we can estimate the minimal and maximal values of the demand taking into account the distribution of clients along the streets. The uncertainty in demand can stem from the uncertainty about the need for or the attractiveness of the facility's service, but can also be due to a variation of the actually present population in the street over the course of the day. For example, in a residential area, the number of inhabitants during the working hours would be smaller than in the evenings.

In the face of this situation of total uncertainty in the demand, we propose to use concepts from robust optimisation to identify suitable facility locations. More concretely, we aim at minimizing the worst-case coverage loss, i.e., minimize the maximal regret. This criterion was first introduced in \citet{Savage51}.  A review of minmax regret optimization as well as a discussion on its potential applications is contained in the book by \citet{kouYu97}. This criterion has been successfully applied to many location problems on networks, see  \citet{AssNorCynAnd17,Averbakh01,AB00,AL04,Conde19,ConLeaPuer18,kouVAiYu93,PueRodTam09}, among others. In case of the MCLP, this criterion is used by \citet{CocSanNor18} for a discrete problem with a finite set of demand points and a finite set of potential locations for facilities. \citet{BerWan11} consider the minmax regret objective for the gradual covering location problem on networks, also assuming that demand only occurs at nodes, but that facilities can be located anywhere along the network. The gradual covering location problem is a generalized version of MCLP where the coverage area is divided in two sectors, one where the demand of nodes within it is completely covered and the other one where the demand is partially covered. An alternative way of considering uncertainty for the discrete MCLP is to use a fuzzy framework, as shown in \citet{AraBlaFer20}.

Our goal in this paper is to solve the minmax regret single-facility MCLP on a network with three main features: edge demand, uncertainty in the intensity of the demand, and a continuous location space. We analyze two cases: i) the realisation of the demand intensity is unknown, but constant on an edge and bounded (from above and below) by two known constant functions; and ii) the demand intensity is given by an unknown linear function on each edge and it is bounded by two known linear functions (see \citet{LopPueRod13,ABL18,KouGan97Chap} for a different way of using uncertain demand intensities at nodes given by linear functions for minmax regret location problems). 

 In the aforementioned examples of locating AED, ATM, bus stops, automated parcel lockers, or bicycle parking racks, the uncertainty in the demand comes from the fact that the number of inhabitatants (potential clients) during the working day is not known. Looking at cities, the type of buildings that one can find along a street are often similar. For example, in the suburbs it is often all single-family homes or low-storey apartment buildings along a street; in more central areas, high-storey apartment buildings or even high rises. If the type of buildings is fairly uniform, we can assume that the lower and upper bounds on the demand are also fairly uniform, i.e., each of them can be modeled by a constant function. 

For such situations, our aim is to find a location such that the worst-case regret for the location is minimal. This assumption makes sense for AED locations, because the health service does not know who is going to need an AED or where, but they are interested in providing this service for all the inhabitants of the city as fast as possible. Similarly, for the location of ATMs, the banks are interested in providing their services to the maximal number of people and avoid losing potential clients.

Regarding linear demand along the edge, this could be the case in the situation of locating bus stops, automated parcel lockers, or bicycle parking racks, where the demand along the streets can vary since either the presence of potential users (tourists, workers, students, etc.) increases when approaching the city center, the business area, or the university campus, or the type of buildings along a street changes gradually, e.g., from single-family homes on one end to row houses on the other end. A suitable way to model such increases could be through a linear function, since, in most cases, these increases are gradual when approaching to those areas. Accordingly, also the lower and upper bounds on the demand intensity could be linear. Maybe not in a strictly linear fashion, e.g., in the case of a change of type of building, but as long as the change  
is not too dramatic, this is a fair assumption. In case of a dramatic change, e.g., from single-family homes to high rises, we would split the edge at the point where the change occurs and treat the two sub-edges separately, reflecting the change also in the lower and upper bound. 
	
Alternative applications can be found if the network represents air ducts of a building (conduits used in heating, ventilation, and air conditioning). For example, one application is to locate an aerosol dispenser in order to disinfect the conduits, see \citet{frazier1990aerosol}. This fits to a MCLP, since aerosols are only effective within a certain distance from their source. Moreover, the demand represents the presence of bacteria, viruses, and fungi along the conduits that should be killed. The actual intensity of demand along the conduits is unknown, but can be estimated by room occupancy rates, uses and types of buildings (healthcare facilities, offices, restaurants, care homes, etc.), and proximity to air inlets and outlets, among others. Since the presence of germs increases when approaching to air inlets and outlets, the demand can be approximated by linear functions.

The remainder of the paper is organised as follows. In Section \ref{sec:ProblemDescription}, we describe the problem.  In Section \ref{sec:ConstantDemand}, we present a polynomial time algorithm to solve the unknown constant demand case and we include a computational study to illustrate the potential and limits of the proposed methodology. In Section \ref{sec:LinearDemand}, we develop a polynomial time algorithm to solve the unknown linear demand case. In both cases, an example illustrating the proposed methodology is shown. 
In Section \ref{sec:Conclusions}, we present our conclusions and an outlook to future research. Finally, there are two appendices. In the first one, we include some known results about the lower envelope of Jordan arcs that will be used to analyse the complexity of the algorithms  proposed in the paper. In the second one, we derive the representation of some functions used in the examples.

\section{Definitions and Problem Description}
\label{sec:ProblemDescription}%

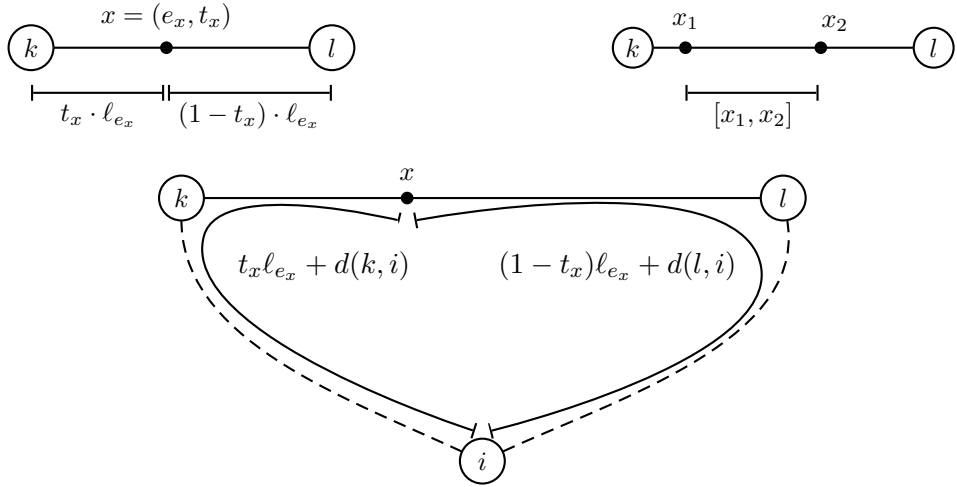
\begin{figure}[htb]
	\centering
	\begin{pspicture}(-0.5,-3.5)(12.5,2.7)
	\psset{radius=0.2, fillstyle=solid}
	
	\Cnode[radius=1.8ex](0,2){i} \rput(0,2){\small $k$}
	\Cnode[radius=1.8ex](4,2){j} \rput(4,2){\small $l$}

	\cnode[fillstyle=solid,fillcolor=black](1.8,2){.07}{sp} \nput{90}{sp}{\small $x=(\ex,t_x)$}
	\psline{|-|}(0,1.4)(1.77,1.4) \rput(0.9,1.1){\small $t_x\cdot {\ell}_{e_x}$}
	\psline{|-|}(4,1.4)(1.82,1.4) \rput(2.9,1.1){\small $(1-t_x)\cdot {\ell}_{e_x}$}
	\ncline{i}{j}

	\Cnode[radius=1.6ex](8,2){k} \rput(8,2){\small $k$}
	\Cnode[radius=1.6ex](12,2){l} \rput(12,2){\small $l$}
	
	\cnode[fillstyle=solid,fillcolor=black](10.5,2){.07}{x} \nput{65}{x}{\small $x_2$}
	\cnode[fillstyle=solid,fillcolor=black](8.7,2){.07}{sp} \nput{90}{sp}{\small $x_1$}
	\ncline{k}{l}
	\psline{|-|}(8.7,1.4)(10.47,1.4) \rput(9.6,1.1){\small $[x_1,x_2]$} 
	
	\Cnode[radius=1.8ex](2,0){i} \rput(2,0){\small $k$}
	\Cnode[radius=1.8ex](10,0){j} \rput(10,0){\small $l$}
	\Cnode[radius=1.8ex](6,-3.5){x} \rput(6,-3.5){\small $i$}
	
	\cnode[fillstyle=solid,fillcolor=black](5,0){.07}{sp} \nput{90}{sp}{\small $x$}
	
	\ncline{i}{j}
	\nccurve[angleA=155, angleB=-90, linestyle=dashed]{x}{i}
	\nccurve[angleA=25, angleB=-80, linestyle=dashed]{x}{j}
	
	\pscurve{|-|}(5.9,-3.1)(2.3,-0.5)(4.9,-0.3) \rput(3.9,-0.9){$t_x {\ell}_{e_x}+d(k,i)$}
	\pscurve{|-|}(6.1,-3.1)(9.7,-0.89)(5.1,-0.3) \rput(7.8,-0.9){$(1-t_x){\ell}_{e_x}+d(l,i)$}
	\end{pspicture}
	\caption{Illustration for a point on an edge, a subedge, and $d(x,i)$.}
	\label{fig:ProbDescr:Example:Notations}
\end{figure}

Let $G=(V,E)$ be an undirected graph with node set $V=\{1, \ldots, n\}$ and edge set $E$, where $|E|=m$.
Every edge $e=[k,l] \in E$, $k,l \in V,$ has a positive length  
$\ell_e$ and is assumed to be rectifiable. $G$ will also denote the continuum set of all points of the network. In addition, in a slight abuse of notation, we use $\ex$ to represent an edge that contains point $x$.
For an edge  $\ex=[k,l] \in E$, a point $x \in \ex$ is also represented as $x = ( \ex,t_x)$, $0 \le t_x \le 1$, where $t_x$ is the relative distance of $x$ from $k$ with respect to ${\ell}_{e_x}$, i.e., the distance from $x$ to node $k$ is $t_x {\ell}_{e_x}$ and to node $l$ is $(1-t_x){\ell}_{e_x}.$  Then, for a node  $i \in V$, the distance from $  x$ to $ i$ is defined as $ d(x,i)  = \min \{ t_x\cdot {\ell}_{e_x} + d(k,i),\, (1-t_x) \cdot {\ell}_{e_x} + d(l,i)   \}$, where $d(k,i)$ and $d(l,i)$ are the lengths of the shortest paths connecting  $k$ with $i$ and $l$ with $i$, respectively.
The distance $d(x,y)$ between two points $x,y\in G$ is defined analogously. 
Let $(e,[t_1,t_2])=[x_1,x_2]=\{x\in e=[k,l]: t_1\ell_e\leq d(x,k)\leq t_2\ell_e, (1-t_2)\ell_e\leq d(x,l)\leq (1-t_1)\ell_e \}$ be a closed subedge of $e$, where $x_1=(e,t_1)$ and $x_2=(e,t_2).$ 
Besides, we are given a fixed coverage radius $R > 0$. See Figure~\ref{fig:ProbDescr:Example:Notations} for an illustration of a point $x$ on $\ex$ (top left-hand side), a subedge $[x_1,x_2]$ (top right-hand side), and of $d(x,i)$ (at the bottom).

Finally, for each edge $e=[k,l] \in E$ we are given two non-negative continuous functions $lb_e: [0,1] \to \Real_0^+$ and $ub_e: [0,1] \to \Real_0^+$ that specify the minimal and maximal demand along the edge.
A specific \emph{demand realisation} or \emph{scenario} on $e$ is denoted by the continuous function $w_e: [0,1] \to \Real_0^+$ where $lb_e(t)\le w_e(t) \le ub_e(t)$, $t \in [0,1]$. For short, we sometimes write $lb_e \le w_e \le ub_e$ for the latter condition.
We define by $lb=(lb_e)_{e \in E}$ the vector of lower bound functions on the network; analogously, we define $ub=(ub_e)_{e \in E}$ and $w=(w_e)_{e \in E}$. Abusing notation, we will abbreviate the corresponding condition for the demand realisations by $lb \le w \le ub$.
\medskip

Let $x \in G$ be a point on the network. We say that a point $z \in G$ is \emph{covered} by a facility at $x$, if $d(x,z) \le R$. Let $C(x) := \{ z \in G \mid d(x,z) \le R \}$ be the \emph{coverage area} of $x$, i.e., the whole set of points of $G$ covered by $x$ and let $C_e(x):= \{ z \in e \mid d(x,z) \le R \}$ be the coverage area of $x$ on $e$, i.e., the set of points of $e$ covered by $x$. 
The total amount of covered demand on an edge $e \in E$ by a facility at $x$ for a specific demand realisation $w$ is given by
\begin{equation}
\label{eq:ProbDescr:CoverageEdge}
g_e(x,w_e) \:=\: \int_{y=(e,t) \in C_e(x)}\, w_e(t)\, dt, 
\end{equation}
and the total amount of covered demand on the entire network by
\begin{equation}
\label{eq:ProbDescr:CoverageTotal}
g(x,w) \:=\: \sum_{e \in E}\, g_e(x,w_e) \,.
\end{equation}
The \emph{worst-case} or \emph{maximal regret} of choosing $x$ over all possible demand realisations is defined as
\begin{equation}
\label{eq:ProbDescr:MaxRegret}
r(x) \::=\: \max_{lb \le w \le ub}\, \left( \max_{y \in G}\, g(y,w) \,-\, g(x,w) \right).
\end{equation}
The realisation $w$ and the point $y$ maximizing the right-hand side  of the expression above are called the \emph{worst-case realisation} and the \emph{worst-case alternative}, respectively, for $x$.
The \emph{Minmax Regret Maximal Covering Location Problem with Edge Demand \name} can now be  formulated as follows: 
\begin{equation}
\label{eq:ProbDescr:MinMaxRegret}
 r^* \::=\: \min_{x \in G}\, \max_{lb \le w \le ub}\, \left( \max_{y \in G}\, g(y,w) \,-\, g(x,w) \right) .
\end{equation}
\medskip

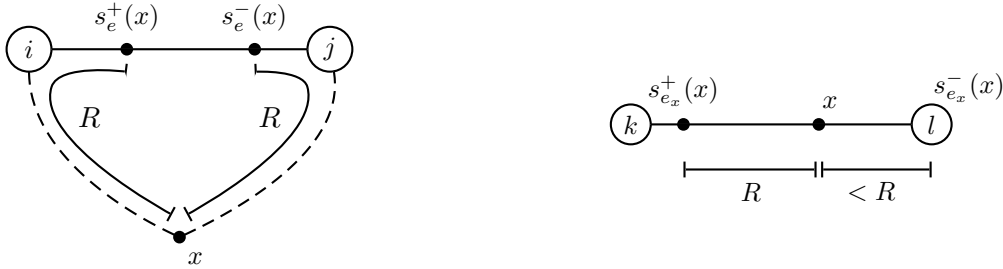
\begin{figure}[htb]
	\centering
	\begin{pspicture}(-0.5,0.2)(12.5,3.7)
	\psset{radius=0.2, fillstyle=solid}
	
	\Cnode[radius=1.8ex](0,3){i} \rput(0,3){\small $i$}
	\Cnode[radius=1.8ex](4,3){j} \rput(4,3){\small $j$}
	
	\cnode[fillstyle=solid,fillcolor=black](2,0.5){.07}{x} \nput{-55}{x}{\small $x$}
	\cnode[fillstyle=solid,fillcolor=black](1.3,3){.07}{sp} \nput{90}{sp}{\small $s_{e}^+(x)$}
	\cnode[fillstyle=solid,fillcolor=black](3,3){.07}{sm} \nput{90}{sm}{\small $s_{e}^-(x)$}
	
	\ncline{i}{j}
	\nccurve[angleA=155, angleB=-90, linestyle=dashed]{x}{i}
	\nccurve[angleA=25, angleB=-80, linestyle=dashed]{x}{j}
	
	\pscurve{|-|}(1.9,0.8)(0.3,2.5)(1.3,2.7) \rput(0.8,2.1){$R$}
	\pscurve{|-|}(2.1,0.8)(3.7,2.5)(3.0,2.7) \rput(3.2,2.1){$R$}
	
	\Cnode[radius=1.6ex](8,2){k} \rput(8,2){\small $k$}
	\Cnode[radius=1.6ex](12,2){l} \rput(12,2){\small $l$}
	
	\cnode[fillstyle=solid,fillcolor=black](10.5,2){.07}{x} \nput{65}{x}{\small $x$}
	\cnode[fillstyle=solid,fillcolor=black](8.7,2){.07}{sp} \nput{90}{sp}{\small $s_{\ex}^+(x)$}
	\rput(12.5,2.5){\small $s_{\ex}^-(x)$}
	\ncline{k}{l}
	\psline{|-|}(8.7,1.4)(10.47,1.4) \rput(9.6,1.1){\small $R$}
	\psline{|-|}(10.53,1.4)(12,1.4) \rput(11.2,1.1){\small $<R$}
	\end{pspicture}
	\caption{Illustration of the edge coverage points $s_e^+(x)$ and $s_e^-(x)$ (\citet{BerKalKra16}).}
	\label{fig:ProbDescr:Example:EdgeCoveragePoints}
\end{figure}

Next, we will recall several definitions and results that were  given in \citet{BerKalKra16} for the case without uncertainty.
Let $x = (\ex,\tx)$ with $\ex = [k,l] \in E$ and $\tx \in [0,1]$. For an arbitrary edge $e = [i,j] \in E$, $e \neq \ex$, let
\begin{align*}
s_e^+(x) & \:=\: \min \left\{1,\, \max \left\{0,\, \frac{R- d(x,i)}{\ell_e} \right\} \right\} \quad \text{and} \\[0.5ex]
s_e^-(x) & \:=\: \max \left\{0,\, \min \left\{1,\, 1-\frac{R- d(x,j)}{\ell_e} \right\} \right\} \,.
\end{align*}
If $s_e^+(x) \in (0,1)$, then $(e,s_e^+(x))$ is the point on $e$ for which the shortest path from this point to $x$ via the end node $i$ has a length of exactly $R$. Moreover, all points on the subedge $(e,[0,s_e^+(x)])$ are covered by $x$. The same interpretation holds for $s_e^-(x)$ with respect to node $j$. See the picture on the left-hand side in Figure~\ref{fig:ProbDescr:Example:EdgeCoveragePoints}.
For $e = \ex$ we define
\begin{align*}
s_{\ex}^+(x) & \:=\: \max \left\{0,\, \frac{d(x,k)-R}{ {\ell}_{e_x}} \right\} \qquad \text{and} \qquad s_{\ex}^-(x) \:=\: \min \left\{1,\, 1-\frac{ d(x,l)-R}{ {\ell}_{e_x}} \right\} \,.
\end{align*}
Hence, all demand points on the subedge $({\ex},[s_{\ex}^+(x),s_{\ex}^-(x)])$ of ${\ex}$ will be covered by a facility at $x$. See the picture on the right-hand side in Figure~\ref{fig:ProbDescr:Example:EdgeCoveragePoints} for an illustration. In the following we call $s_e^+(\cdot)$ and $s_e^-(\cdot)$ the \emph{edge coverage functions}. 

In the following lemma, we sum up the calculation of the coverage area of $x$:

\begin{lemma}[\citet{BerKalKra16}]
\label{lem:ProbDescr:CoverageArea}
For $e \in E \setminus \{ \ex \}$, we have
\[ C_e(x) \:=\: \begin{cases} e, & \text{ if } s_e^-(x) \le s_e^+(x), \\
(e,[0,s_e^+(x)] \cup [s_e^-(x),1]), & \text{ if } 0 < s_e^+(x) < s_e^-(x) < 1, \\
(e,[0,s_e^+(x)]), & \text{ if } 0 < s_e^+(x) < 1 = s_e^-(x), \\ 
(e,[s_e^-(x),1]), & \text{ if } s_e^+(x)= 0 < s_e^-(x) < 1, \\
\emptyset, & \text{ if } s_e^+(x)=0,\, s_e^-(x)=1,
\end{cases}
\]
and $C_{\ex}(x) = (\ex,[s_{\ex}^+(x),s_{\ex}^-(x)])$.
\end{lemma}

Using this, we classify all edges $E \setminus \{\ex\}$ as covered, uncovered, and partially covered as follows:
\begin{align*}
E^c(x) & \:=\: \{ e \in E \setminus \{\ex\} \mid s_e^-(x) \le s_e^+(x) \}, \\[0.5ex]
E^u(x) & \:=\: \{ e \in E \setminus \{\ex\} \mid s_e^+(x) = 0 \text{ and } s_e^-(x) = 1 \}, \\[0.5ex]
E^p(x) & \:=\: E \setminus \{E^c(x) \cup E^u(x) \cup \{\ex\}\}.
\end{align*}
Edge $\ex$ will either be fully or partially covered.
The total coverage can now be written as
\begin{align}
\label{eq:ProbDescr:CoverageTotalReformulation}
g(x,w) \:=\: & \int_{s_{\ex}^+(x)}^{s_{\ex}^-(x)}\, w_{\ex}(u)\, du \,+\, \sum_{e \in E^c(x)}\, \int_0^1\, w_e(u)\, du \nonumber \\
& \,+\, \sum_{e \in E^p(x)}\, \left( \int_0^{s_e^+(x)} w_e(u)\, du + \int_{s_e^-(x)}^1 w_e(u)\, du \right)
\,.
\end{align}

	\subsection{Singularity points}
	 In the current subsection, we analyse the singularity points on a network  that will be useful to analyse the behaviour of the objective function for the problem under study.
 A point $ x = (\ex,\tx)$, $\ex= [k,l] \in E$ is called a bottleneck point of node  $i$, if $\tx\cdot  {\ell}_{e_x} +d(k,i) \,=\, (1-\tx) \cdot {\ell}_{e_x} + d(l,i) $ with $0 < \tx < 1$.
	We denote by $BP_{\ex}$ and $BP$ the set of all bottleneck points on edge $\ex$ and all bottleneck points of the network, respectively. 
	
	In the following, we denote by $NP_{\ex} = \{ x \in \ex \mid \exists\,  i \in V:\, d(x,i)=R \}$ the set of all network intersect points (NIPs) on an edge $\ex$ and by $NP = \bigcup_{\ex \in E}\, NP_{\ex}$ the set of all NIPs on the network (\citet{ChuMea79}).
	
 The bottleneck points and  network intersect points determine the breakpoints of the edge coverage functions, as stated in the following result. 
		\begin{lemma}[\citet{BerKalKra16}]
			\label{lm:breakpoints}
			The edge coverage functions $s_e^+(x)$ and $s_e^-(x)$, $e \in E$, are continuous and piecewise linear functions over $x\in \ex$ with a constant number of pieces.  Breakpoints of these two functions are either bottleneck points or network intersect points (associated with some of the two endnodes of $e$).
		\end{lemma}
	
	A point $x \in G$ is called \emph{exact coverage point} if $0 < s_e^+( x) = s_e^-( x) < 1$ for some $e \in E$ 
	and there exists $\epsilon_1,\epsilon_2\geq 0$, $\epsilon_1+\epsilon_2> 0$, such that all points $x^\prime \in (x-\epsilon_1, x) \cup (x,x+\epsilon_2),$ no longer completely cover $e$ (\citet{BerKalKra16}). We denote by $EP_{\ex}$ the set of all exact coverage points on an edge $\ex$ and by $EP = \bigcup_{\ex \in E}\, EP_{\ex}$ the set of all exact coverage points on the network.
	
	Finally, for $\ex= [k,l]\in E$, we define the set of \emph{partition points}
	\[ PP_{\ex} \::=\: \{ k,l\} \cup NP_{\ex} \cup BP_{\ex} \cup EP_{\ex} \qquad \text{and} \qquad PP \::=\: \bigcup_{\ex \in E}\, PP_{\ex} \,.
	\]
	For two consecutive partition points $z^1$, $z^2 \in PP_{\ex}$, we call the subedge $[z^1,z^2] \subseteq \ex$ a \emph{partition edge}.
	Concerning the cardinality of $PP$,  let $e,\ex \in E$, such that $e\neq \ex$. Each node of the network induces at most one bottleneck point and at most two network intersect points on $\ex$. Similarly, we compute the number of exact coverage points. By Lemma \ref{lm:breakpoints},  $s_e^+(x)$ and $s_e^-(x)$, $e \in E$, are continuous and piecewise linear functions over $x\in \ex$ with a constant number of pieces. Hence, there is a constant number of  solutions on $\ex$ of $s_e^+(x) = s_e^-( x)$, i.e.,  there is a constant number of exact coverage points. Therefore, there are at most $\mathcal{O}(m)$ partition points  on each $\ex\in E$ and $\mathcal{O}(m^2)$ on the whole network.
With the above definitions, we can formulate the following result:
\begin{prop}[\citet{BerKalKra16}]
\label{prop:ProbDescr:PartitionPoints}
Let $\ex \in E$ and $z^1,z^2 \in PP_{\ex}$ be two consecutive partition points on $\ex$. Moreover, let $w$ be a non-negative continuous demand function.
\begin{enumerate}
\item
The sets $E^c(x)$, $E^u(x)$, and $E^p(x)$ are identical for all $x \in [z^1,z^2]$.

\item
The edge coverage functions $s_e^+(x)$ and $s_e^-(x)$ have a unique linear representation in $x$ over $[z^1,z^2]$ for each partially covered edge $e \in E^p(x)$.

\item
All functions $g_e(x,w_e)$, $e \in E$, have a unique representation in $x$ over $[z^1,z^2]$.
\end{enumerate}
\end{prop}

\noindent
 A summary of the notation introduced in this section can be seen in Table \ref{tab:notation}.
\begin{table}[htb]
        \caption{Notation used in the paper.}
    \centering
     \resizebox{\hsize}{!}{
    \begin{tabular}{ll}
    \hline
     $BP_e, BP$ & the set of all bottleneck points on edge $e$ and on $G$, respectively.\\
      $c_e(x)$ &the parts per unit of coverage function of $e.$\\
         $C_e(x),C(x)$ & the coverage area of $x$ on $e$ and on $G$, respectively.\\
         $d(x,y)$ & distance between two points $x,y\in G.$\\
         $e=[k,l]$ & edge, where $k,l\in V.$ \\
          $\ex$ & edge that contains point $x.$ \\
         $E$  & the set of edges.\\
         $E^c(x),E^u(x),E^p(x)$& the set of edges that are covered, uncovered and partially covered by $x.$\\
         $EP_e, EP$ & the set of all exact coverage points on edge $e$ and on $G$, respectively.\\
         $G=(V,E)$ & network with node set $V$ and edge set $E$ and the continuum set of points of it.\\ 
         $IC_e, IC$ & the set of all identical coverage points on edge $e$ and on $G$, respectively.\\
         $\ell_{e}$& length of edge $e\in E.$\\
         $lb_e, ub_e$& the functions that specify the minimal and maximal demand over edge $e.$\\
         $m$&number of edges.\\
         $n$&number of nodes.\\
         $NP_e, NP$ & the set of all network intersect points on edge $e$ and on $G$, respectively.\\
         $PP_e, PP$ & the set of all partition points on edge $e$ and on $G$, respectively.\\
         $r(x)$ & the maximal regret of location $x$ over all possible demand realisations.\\
         $r_e(x,y), r(x,y)$ & the maximal regret for $x$ with respect to $y$ on an edge $e$ and on $G$, respectively.\\
         $R$ & the coverage radius.\\
         $s_e^{+}(x),s_e^{-}(x)$ &the edge coverage functions of $e.$\\
         $x=(e_x,t_x)$& point in edge $e_x=[k,l]$ where $t_x$ is the relative distance of $x$ from $k$ with respect to $l_{\ex}.$\\
         $[x_1,x_2]=(e,[t_1,t_2])$ & subedge of $e,$ where $x_1=(e,t_1),$ $x_2=(e,t_2),$ and $0\leq t_1<t_2\leq 1.$\\
        $V$ & the set of nodes. \\
         $w_e$ & demand realization on edge $e.$\\
         \hline
    \end{tabular}
    }

    \label{tab:notation}
\end{table}


\section{The Unknown Constant Demand Case}
\label{sec:ConstantDemand}%

In this section, we solve the single facility \name with unknown constant demand realisations $w_e(t)$, i.e., we assume that the demand realisation along the edges is an unknown but constant function bounded by known lower and upper constant functions.  We will derive theoretical properties of the solution with the objective of providing an exact algorithm to solve the problem in polynomial time.

We first look at the covered demand. 
Slightly abusing notation, we denote by $ub_e$, $lb_e$, and $w_e$ the constant value of the lower and upper bound function and the demand realisation, respectively, on $e$.
Let $x = (\ex,\tx)$ be given with $\ex = [k,l] \in E$ and $\tx \in [0,1]$. Using \eqref{eq:ProbDescr:CoverageTotalReformulation}, we obtain
\begin{align}
g(x,w) & \:=\: w_{\ex} (s_{\ex}^-(x) - s_{\ex}^+(x)) \,+ \sum_{e \in E^c(x)}\, w_e \,+\, \sum_{e \in E^p(x)}\, w_e \left( 1 - (s_e^-(x) - s_e^+(x)) \right).
\end{align}

\noindent
By Proposition~\ref{prop:ProbDescr:PartitionPoints}, we obtain that the function $g(x,w)$ is linear in $x$ over each  partition edge of $\ex$ and it allows us to get  the following result.

\begin{prop}[\citet{BerKalKra16}]
\label{prop:Constant:FDS}
$PP$ is a finite dominating set for the single facility maximal covering location problem with constant demand on  each edge.
\end{prop}

\begin{observation}
\label{obs:Constant:DemandIndependentFDS}
The finite dominating set does not depend on the actual realisation $w$.
\end{observation}

\noindent
From Proposition~\ref{prop:Constant:FDS} and  Observation~\ref{obs:Constant:DemandIndependentFDS}, we immediately obtain

\begin{lemma}
\label{lem:Constant:Regret}
The maximal regret of choosing $x$ over all possible constant demand realisations  on each edge is $r(x) \:=\: \displaystyle\max_{lb \le w \le ub}\, \left(\max_{y \in PP}\, g(y,w) \,-\, g(x,w)\right).$ 
\end{lemma}

\noindent
In the following, we denote
\[ c_e(x) \::=\: \begin{cases} 1, & \text{ if } e \in E^c(x), \\
1 - (s_e^-(x) - s_e^+(x)), & \text{ if } e \in E^p(x), \\
0, & \text{ if } e \in E^u(x), \\
s_{\ex}^-(x) - s_{\ex}^+(x), & \text{ if } e = \ex,
\end{cases}
\]
as the \emph{parts per unit of coverage} of the respective edge.
Observe that $0 \le c_e(x) \le 1$ and $g(x,w)=\sum_{e \in E} w_e\,c_e(x)$.
 For a given partition point $y\in PP,$ the maximal regret for $x$ with respect to $y$ on an edge $e\in E$ can be computed as:
\begin{equation}
r_e(x,y)\::=\: \max_{lb_e \le w_e \le ub_e}\, g_e(y,w_e)-g_e(x,w_e).
\end{equation}
Therefore, the maximal regret for $x$ with respect to $y$ on the entire network is given by:
\begin{align}
\label{eq:Constant:MaxRegretFixedY}
r(x,y) & \::=\:\max_{lb \le w \le ub}\, (g(y,w)-g(x,w))
 \:=\: \max_{lb \le w \le ub}\, \sum_{e \in E}(g_e(y,w_e)-g_e(x,w_e))  \nonumber \\ &  \:=\:\sum_{e \in E}\max_{lb_e \le w_e \le ub_e}\,(g_e(y,w_e)-g_e(x,w_e)) \:=\:\sum_{e \in E} r_e(x,y) \nonumber
\\&\:=\: \sum_{e \in E}\,\max_{lb_e \le w_e \le ub_e}\,  w_e\,(c_e(y)-c_e(x)) \nonumber \\
  & \:=\: \sum_{e \in E: c_e(y) \ge c_e(x)}\, ub_e\,(c_e(y)-c_e(x)) \,+\, \sum_{e \in E: c_e(y) < c_e(x)}\, lb_e\,(c_e(y)-c_e(x)).
\end{align}
\begin{theorem}
\label{theo:Constant:MaxRegret}
The maximal regret of $x \in G$ is given by
\begin{align}
\label{eq:Constant:MaxRegret}
r(x) \:=\: \max_{y \in PP}\,  \sum_{e \in E}\left(g_e(y,w_e^{x,y}) - g_e(x,w_e^{x,y})\right)\,,
\end{align}
where
\[ w_e^{x,y} \:=\: \begin{cases}
                     ub_e, & \text{if } c_e(y) \ge c_e(x), \\
                     lb_e, & \text{otherwise,}
                   \end{cases}
\]
is the worst-case demand realisation with respect to $x$ and $y$.
\end{theorem}

\ProofNoNL
Let $y^*$ and $w^*$ be the worst-case alternative and worst-case demand realisation, respectively, for $x\in G,$ see expression~\eqref{eq:ProbDescr:MaxRegret}, i.e., $r(x) = g(y^*,w^*) \,-\, g(x,w^*)$. From Proposition~\ref{prop:Constant:FDS} and  \eqref{eq:Constant:MaxRegretFixedY} we can assume without loss of generality that $y^* \in PP$ and $w^*=w^{x,y^*}$, respectively.
Thus,
\begin{align*}
r(x) & \:=\: \max_{y \in PP}\, \left(g(y,w^{x,y}) - g(x,w^{x,y})\right) \,.
\end{align*}\EndProofNoNL

Having established an easy way to compute the maximal regret of $x$, the only open question is how to optimise $x$ over $\ex$.
Let $[z^1,z^2]$ be a partition edge of $\ex$ and $x \in \text{int}([z^1,z^2])$. Moreover, let $y \in PP$ be given.
From Proposition~\ref{prop:ProbDescr:PartitionPoints} we know that the edge coverage functions are linear over $[z^1,z^2]$ and that the sets $E^p(x)$, $E^c(x)$, and $E^u(x)$ are identical for all $x \in [z^1,z^2]$. What, however, can change over $[z^1,z^2]$ is the demand realisation $w^{x,y}$ that yields the maximal regret for $x$ and $y$.
In that case, there must exist an edge $e\in E$ and a point $z \in [z^1,z^2]$ such that $z$ provides the exact same coverage for $e$ as $y$, i.e., $ c_e(z) \:=\: c_e(y).$ To see that, we next look at the four possible cases (in the remaining cases, i.e., $e\in E^c(z)\cap E^c(y)$ and $e\in E^u(z)\cap E^u(y)$, $r_e(z,y)=0$ for any $z\in [z^1,z^2]$):

\noindent If $e\in E^p(z) \cap E^p(y),$ or $z\in e$ and $y\in e$, then:  \[c_e(z)  \:=\: c_e(y)\quad \Leftrightarrow \quad s_e^-(z) - s_e^+(z) \:=\: s_e^-(y) - s_e^+(y).\] If $e\in E^p(z)$ and $y\in e,$ then:  \[c_e(z)  \:=\: c_e(y) \quad\Leftrightarrow \quad 1-(s_e^-(z) - s_e^+(z))  \:=\: s_e^-(y) - s_e^+(y).\] Finally, if $z\in e$ and $e\in E^p(y),$ then:  \[c_e(z)  \:=\: c_e(y) \quad\Leftrightarrow \quad s_e^-(z) - s_e^+(z)  \:=\: 1-(s_e^-(y) - s_e^+(y)).\]

As the edge coverage functions, $s_e^+(z)$ and $s_e^-(z),$ for $z\in[z^1,z^2]$ are linear functions, the equation $c_e(z)=c_e(y)$ will have either  at most one solution  or a subinterval contained in $[z^1,z^2]$ as solution. In the former case, if a solution exists it will be called \textit{identical coverage point} with respect to $y$ and $e$ and in the later case,  \wloge, the extremes of this subinterval will be called the identical coverage points of $[z^1,z^2]$ with respect to $y$ and $e$. Furthermore, in both cases the solution is independent of $w$.
We define $IC_{\ex}(y)$ as the set of all identical coverage points of $G$ induced by $y$ and $\ex$. Moreover, we define $IC = \bigcup_{\ex \in E, y \in PP}\, IC_{\ex}(y)$.

\begin{theorem}
\label{theo:Constant:Linearity}
Let $z^1$, $z^2$ be two consecutive elements of $PP \cup IC$ on an edge $\ex \in E$ and $y \in PP$. Then
\begin{enumerate}
\item
$g(y,w^{x,y}) - g(x,w^{x,y})$ is linear in $x$ over $[z^1,z^2]$.
\item
$r(x)$ is piecewise linear and convex in $x$ over $[z^1,z^2]$.
\end{enumerate}
\end{theorem}

\ProofNoNL
Let $z^1$, $z^2$ be two consecutive elements of $PP \cup IC$ and let $x' \in \text{int}([z^1,z^2])$. Moreover, let $y \in PP$ and $w^{x',y}$ the corresponding worst-case demand realisation. By definition, we can assume without loss of generality that $w^{x',y}=w^{x,y}$ for all $x \in [z^1,z^2].$ 
From Proposition~\ref{prop:ProbDescr:PartitionPoints} we know that the edge coverage functions are linear over $[z^1,z^2]$ and that the sets $E^p(x)$, $E^c(x)$, and $E^u(x)$ are identical for all $x \in [z^1,z^2]$. Hence, $g(y,w^{x,y}) - g(x,w^{x,y})$ is linear in $x$ over $[z^1,z^2]$. The second result then follows immediately from  the previous statement and  Theorem~\ref{theo:Constant:MaxRegret}.
\EndProofNoNL
	
\begin{observation}
	Theorem~\ref{theo:Constant:Linearity} gives us a partition in the domain of the function $r(x)$ where it is convex. The optimal solution of the problem is the minimum of this function, which can be computed as the minimum of the minima in each subedge. However, we propose an alternative strategy for computing the optimal solution of this problem whose complexity is smaller than that of this mentioned procedure.
\end{observation}	
	\begin{theorem}\label{tm:Constant:rxy}
		For a given $e\in E$ and $y\in PP,$  $r_e(x,y)$ is a piecewise linear function for $x\in \ex$  with a constant number of pieces and $r(x,y)$ is a piecewise linear function with $\mathcal{O}(m)$ number of pieces for $x\in \ex.$
	\end{theorem}
	\ProofNoNL
	Let $\ex, e \in E$ and $y \in PP$. By  Lemma~\ref{lm:breakpoints}, we know that the edge coverage functions $s_{e}^+(x),$ and $s_{e}^-(x)$ 
	are piecewise linear functions over $\ex \in E$ with a constant number of breakpoints for $x\in\ex$. 
	Hence, by definition, for each $e\in E$, $c_e(x)$ is also a piecewise linear function with a constant number of breakpoints on $\ex$. Thus using \eqref{eq:Constant:MaxRegretFixedY} and the fact that the number of identical coverage points on $\ex$ for an edge $e\in E$ is constant, by Theorem~\ref{theo:Constant:MaxRegret} there is a constant number of possible worst-case demand realisation, $w_e^{x,y}$ and then,  $r_e(x,y)$ is a piecewise linear function with a constant number of pieces.
	\EndProofNoNL

\noindent
Theorems~\ref{theo:Constant:MaxRegret} and~	\ref{tm:Constant:rxy} give rise to Algorithm~\ref{Al:const2} to find an optimal solution for the \name  with unknown constant demand realisations. 
	
	\begin{algorithm2e}[h]\label{Al:const2}
		\DontPrintSemicolon \SetAlFnt{\small\sl}
		\SetAlCapFnt{\small\sl} \AlCapFnt
		
		\caption{Optimal algorithm for the single facility \name with unknown constant demand}
		\label{algo:LinearConst}
		\BlankLine 
		
		\KwIn{Network $G=(V,E)$; lower and upper bounds $lb_e$, $ub_e$, respectively for $e\in E$; coverage radius $R > 0$.}
		\BlankLine 
		\KwOut{Optimal solution $x^*$.}
		\BlankLine \BlankLine
		
		\nl Determine the set $PP$ of partition points. \;
		\BlankLine
		
		\BlankLine
		\nl \ForEach{$\ex\in E$}{
			\nl \ForEach{$y \in PP$}{
	
				\nl Derive the representation of $r_e(x,y),$ for each $e\in E, x\in \ex$.
				\BlankLine 
				\nl Obtain the representation of $r(x,y),$ for $x\in \ex$.
				\BlankLine
			}	
			
			\BlankLine
			\nl Obtain the upper envelope of $r(x,y)$ for $y\in PP$, i.e., obtain $r(x)$ for $x\in\ex$.
			\BlankLine
			\nl Find the minimum $x_{\ex}^{\prime}$ of $r(x)$ over $\ex$.\;
			\BlankLine
			
			\nl \lIf{$r(x_{\ex}^{\prime})< r(x^*)$}{set $x^* := x_{\ex}^{\prime}$, $r(x^*)= r(x_{\ex}^{\prime})$.}

		}
		\BlankLine
		
		\nl \Return\ $x^*$.\;
	\end{algorithm2e}
	\medskip
 For computing the complexity of the Algorithm~\ref{Al:const2} some technical results are needed about the complexity of computing an upper envelope of Jordan arcs. This complexity is expressed in terms of $\lambda_s(n)$, the maximum length of a Davenport-Schinzel sequence of order $s$ on $n$ symbols, see \citet{vsarir1995davenport}. The results used in the following proof can be found in the Appendix.  
	
\begin{theorem}
	The single facility \name with unknown constant demand realisations on each edge can be solved exactly in $\mathcal{O}(m\lambda_3(m^3))$ time using Algorithm~\ref{Al:const2}.
\end{theorem}

	\ProofNoNL 
	As the time to compute each partition point is constant, Step~1 requires $\mathcal{O}(m^2)$ time (\citet{BerKalKra16}).
	
	Since there are $m$ edges, $r(x,y)$ has $m$ addends,  each one being a piecewise linear function with a constant number of pieces for each $y\in PP$ and $e\in E,$ (Theorem~\ref{tm:Constant:rxy}).   Step~4 takes constant time for each $e\in E$ and fixed $y\in PP$. Besides Step~5  generates a function $r(x,y)$ piecewise linear with $\mathcal{O}(m)$ breakpoints.

	Step~6 obtains the upper envelope of  $\mathcal{O}(m^3)$ line segments, i.e., $\mathcal{O}(m)$ line segments for each $y\in PP$. Since they are line segments, we have $s=1$. Therefore, this step takes $\mathcal{O}(\lambda_{2}(m^3)\log m) =\mathcal{O}(m^3\log m)$ time (Theorem~\ref{prop:compcostupperenvelope}). In addition, the complexity of computing the minimum in Step~7 is dominated by the complexity of the upper envelope. Hence, this step takes  $\mathcal{O}(\lambda_{3}(m^3))$ time and Step~8 takes constant time.
	
	 Step~4 is executed once for each $\ex\in E, y\in PP$ and $e\in E$, i.e., $\mathcal{O}(m^4)$ times. Moreover Step~5 is executed once for each $\ex\in E,y\in PP$, i.e.,  $\mathcal{O}(m^3)$.  Similarly, Step~6-8 are executed once for each  $\ex\in E$, i.e., $\mathcal{O}(m)$ times. Thus, the overall complexity of the algorithm is $\mathcal{O}(m\lambda_{3}(m^3)).$ 
\EndProofNoNL
 
\begin{observation}
	Recall that $\lambda_{3}(m)=\Theta(m\alpha(m))$ where  $\alpha(m)$ is the functional inverse of the Ackermann's function. A weaker but simpler upper bound is $\lambda_{3}(m)=\mathcal{O}(m\log m),$ see \citet{vsarir1995davenport} for further details.
\end{observation}

 \begin{example}\label{example_cd} 
Consider the network depicted in Figure \ref{fig:ExLD1}. For each edge $e\in E,$ its length is printed next to the edge. Let $R=1,$ $lb_{[1,2]}= 3,$ $ub_{[1,2]}= 15,$ $lb_{[2,3]}=1,$ $ub_{[2,3]}= 7$,  $lb_{[1,3]}= 2$, and $ub_{[1,3]}= 8$. The set of partition points is given by $PP=V\cup \{([1,3],1/3), ([1,3],2/3),([2,3],1/2)\}.$ The three partition points not included in $V$, indicated as dots in the figure, are the  bottleneck point $([1,3],2/3)$ and the three network intersect points $(([1,3],1/3),([1,3],2/3)$, and $([2,3],1/2)$.
	
\begin{figure}[htb]
	\centering
	\begin{pspicture}(-0.5,0.2)(12.5,3.7)
	\psset{radius=0.2, fillstyle=solid}
		
	\Cnode[radius=1.8ex](4,3){1} \rput(4,3){\small $1$}
	\Cnode[radius=1.8ex](8,3){2} \rput(8,3){\small $2$}
	\Cnode[radius=1.8ex](6,0.5){3} \rput(6,0.5){\small $3$}
	\cnode[fillstyle=solid,fillcolor=black](4.66,2.16){.07}{pp1}
	\cnode[fillstyle=solid,fillcolor=black](5.33,1.33){.07}{pp2}
	\cnode[fillstyle=solid,fillcolor=black](7.0,1.75){.07}{pp3}
		
	\rput(6,3.3){1}
	\rput(7.35,1.65){2}	
	\rput(4.65,1.65){3}		
	\ncline{1}{2}
	\ncline{1}{3}
	\ncline{3}{2}
		
	\end{pspicture}
	\caption{Network with edge lengths.}
	\label{fig:ExLD1}
\end{figure}
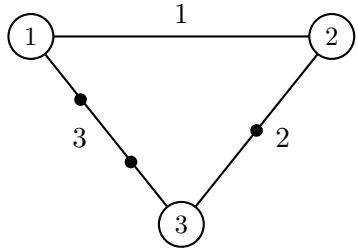

First, we derive the representation of the edge coverage functions and the parts per unit coverage functions. Their expressions can be found in Appendix \ref{sec:app2}. Next, for each edge $e_x,$ we compute the maximal regret for $x\in e_x$ with respect to $y\in PP$: 
$$r(x,y)= \sum_{e \in E: c_e(y) \ge c_e(x)}\, ub_e\,(c_e(y)-c_e(x)) \,+\, \sum_{e \in E: c_e(y) < c_e(x)}\, lb_e\,(c_e(y)-c_e(x))$$ 
In Figure \ref{fig:CDr12}, functions $r(x_1,y),$ for $(x_1,y)\in[1,2]\times PP,$ are depicted. The upper envelope, $r(x),$ is represented as a dotted line. Note that $r(x_1,y_3),$ is equal to $r(x_1,y_4),$ for $y_3=([1,3], 1/3)$ and $y_4=([1,3], 2/3).$ The  reason of it is that the parts per unit of coverage functions are equal. 	As can be observed, the minimum is found in $x^\prime_{[1,2]}=\left([1,2],\frac{2}{3}\right)$ and the regret is $r(x^\prime_{[1,2]})=\frac{13}{9}.$
	
\begin{figure}[h]
	\centering
	\includegraphics[height=0.35\linewidth]{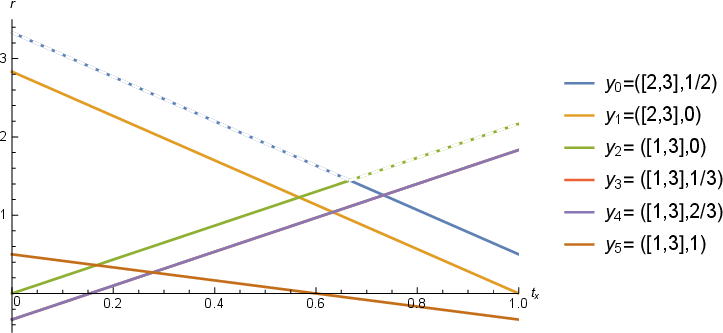}
	\caption{Representation of $r(x_1,y), \; (x_1,y)\in[1,2]\times PP$.}
	\label{fig:CDr12}
\end{figure}
	
Similarly, functions $r(x_2,y),$ for $(x_2,y)\in[2,3]\times PP,$ are depicted in Figure \ref{fig:CDr23}. The upper envelope, $r(x),$ is represented as a dotted line. As before, observe that $r(x_2,y_3),$ is equal to $r(x_2,y_4),$ for $y_3=([1,3], 1/3)$ and $y_4=([1,3], 2/3).$ 	Over this edge, $[2,3]$, the minimum is found in $x^\prime_{[2,3]}=\left([2,3],0\right)$ and the regret is $r(x^\prime_{[2,3]})=\frac{13}{6}.$
		
\begin{figure}[h]
	\centering
	\includegraphics[height=0.35\linewidth]{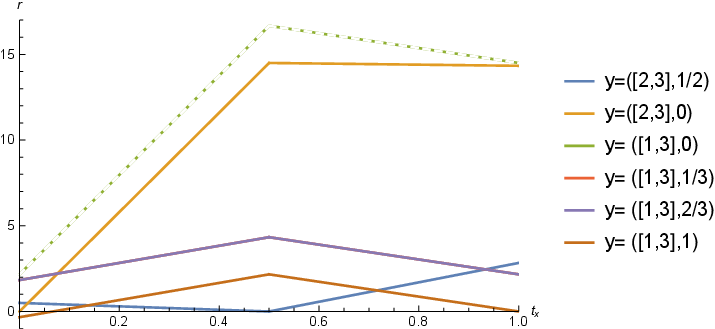}
	\caption{Representation of $r(x_2,y), \; (x_2,y)\in[2,3]\times PP$.}
	\label{fig:CDr23}
\end{figure}

Finally, functions $r(x_3,y),$ for $(x_3,y)\in[1,3]\times PP,$ are represented in Figure \ref{fig:CDr13}. The upper envelope, $r(x),$ is depicted as a dotted line. As before, note that $r(x_3,y_3),$ is equal to $r(x_3,y_4),$ for $y_3=([1,3], 1/3)$ and $y_4=([1,3], 2/3).$ 	Over this edge, $[1,3]$, the minimum is $x^\prime_{[1,3]}=\left([1,3],0\right)$ and the regret is $r(x^\prime_{[1,3]})=\frac{10}{3}.$
	
\begin{figure}[h]
	\centering
	\includegraphics[height=0.35\linewidth]{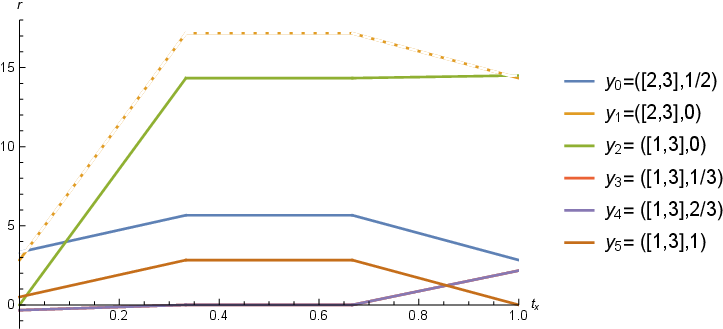}
	\caption{Representation of $r(x_3,y), \; (x_3,y)\in[1,3]\times PP$.}
	\label{fig:CDr13}
\end{figure}

Therefore, the optimal location of the facility with the objective of minimising the maximal regret is $x^\star=\left([1,2],\frac{2}{3}\right)$ and the regret is $r(x^\star)=\frac{13}{9}.$

To highlight the usefulness of our max-regret approach, we also compute the optimal location for a deterministic version of the problem. To that end, we assume that the demand is known and equal to the mean of the upper and lower bound functions over each edge, i.e., $w_e=0.5(ub_e+lb_e)$. Afterwards, we determined the optimal solution $x_D^*$ for this deterministic version using the algorithm in \citet{BerKalKra16}. Finally, we computed the maximal regret of this solution, i.e., $r(x_D^*)$. The optimal solution of the deterministic problem is vertex 2, whose maximal regret is $r(x_D^*)=\frac{13}{6}$. Moreover, we also calculated the amount of covered demand in the deterministic version of the problem for the solution $x_D^*$ as well as for the optimal solution $x^\star$ of the max-regret problem. The covered demand (the objective value of the deterministic model) is 11 and 10.8889, respectively. Therefore, while the amount of covered demand for both locations is almost identical, the maximal regret of the minmax-regret solution is significantly lower than the one for the deterministic solution. The solutions of both models, their maximal regret, and the covered demand (assuming that the demand is deterministic) are given in Table~\ref{tab:summary}.

\begin{table}[htbp]
\caption{Solutions for Example~\ref{example_cd}.}
\centering
{\tablinesep=2ex\tabcolsep=10pt \begin{tabular}{l|r|r|r}
Model & Optimal solution & Maximal Regret & \textit{Covered Demand}  \\
\hline
  \name     &  $\left([1,2],\frac{2}{3}\right)$     &   $\dfrac{13}{9}$    & 10.8889 \\
\hline
   \citet{BerKalKra16}   &  $\left([1,2],1\right)$    &  $\dfrac{13}{6}$     & 11
\end{tabular}%
}
\label{tab:summary}%
\end{table}%
\end{example}

\bigskip
\subsection{Computational experiments}

In this section we present computational experiments for Algorithm~\ref{Al:const2} for randomly generated as well as real-world data sets. Starting with the former, for a given number of vertices $n \in \{40,60,80,100\}$ and a given edge density $p \in \{0.1,0.2,0.3\}$ (density as the percentage of the edges of a complete graph), we randomly generated five connected graphs (in total 60 graphs). For each graph, we first drew the edge lengths from a uniform distribution over $[1,20]$. Afterwards, we computed the distance matrix and checked that the graph did not violate the triangle inequality (if this happened, we replaced the corresponding edge lengths with the shortest path distances and recalculated the distance matrix; we repeated this process until the triangle inequality was satisfied for all edges). In the next step, we randomly generated the values $lb_e$ and $ub_e$ for each edge. To that end, we fixed a value $UB \in \{10,50,100\}$ and drew $lb_e$ and $ub_e$ from a uniform distribution over $[0,UB/2]$ and $[UB/2,UB]$, respectively. Finally, the coverage radius $R \in \{0.1\cdot d^{max}, 0.2\cdot d^{max}, 0.3\cdot d^{max}\}$ is taken as a fixed percentage of the diameter $d^{max} = \max_{i,j\in V} d(i,j)$ of the graph. As a result, for a given combination of $n$ and $p$, i.e., a given graph, we have nine different instances.
The two real-world data sets are based on different street graphs from a German city, with the edge lengths being the length of the street in meters. As we did not have access to actual demand data, we generated upper and lower bounds in the same way as for the randomly generated graphs. Analogously, we chose the coverage radius $R \in \{0.1\cdot d^{max}, 0.2\cdot d^{max}, 0.3\cdot d^{max}\}$.

We compare our algorithm with two alternative versions. The first alternative is the node-restricted (MMR-EMCLP), i.e., we replace $\min_{\in G}$ by $\min_{x \in V}$ in~\eqref{eq:ProbDescr:MinMaxRegret}. The motivation for the second alternative was to get an idea about the value of the stochastic solution for this setting. To that end, we first replaced the unknown constant demand on an edge $e$ by a constant known demand value $w_e = (ub_e+lb_e)/2$. Afterwards, we determined the optimal solution for this deterministic problem using the algorithm in \cite{BerKalKra16}. Finally, we computed the maximal regret of this solution. In the following, we call these two variants simply the \emph{node-restricted} and the \emph{deterministic} algorithm, and our method the \emph{max-regret} algorithm. The corresponding problems and optimal solutions are denoted analogously.
All algorithms were implemented in C++ and run on a Windows 10 Laptop with a i5-8350U CPU with 1.9 GHz and 8 GB RAM.

\subsubsection{Random graphs}

We start with the randomly generated graphs. We first analyse the effect of varying the edge densities. The results are given in Table~\ref{tab:compexp:rand:density}. The first two columns specify the number of nodes and the edge density. The next column shows the average total run time of the max-regret algorithm. The following three columns are for the node-restricted algorithm. The first shows the average relative percentage deviation of the maximal regret of the node-restricted solution with respect to the maximal regret of max-regret solution. The next two columns present the maximum relative percentage deviation and the total run time. The following three columns show exactly the same information, but this time for the deterministic algorithm. 
All averages are taken over the 45 different instances that have the same edge density and number of nodes (five different graphs per $n$ and $p$, and nine different instances per graph).

\begin{table}[htb]
\centering
\small
\caption{Comparison for different edge densities for the random graphs.}
\begin{tabular}{cc|c|ccc|ccc}
& \multicolumn{1}{c}{} & \multicolumn{1}{c}{Max regret} & \multicolumn{3}{c}{Node-restricted} & \multicolumn{3}{c}{Deterministic} \\
Nodes & Density &   Time  &  Avg Dev & Max Dev & Time  &  Avg Dev & Max Dev & Time \\ \hline
40 & 0.1 &     2.10 &   0.8\% &  7.6\% & 0.05 &   21.8\% & 318.2\% & 0.03  \\
40 & 0.2 &    16.77 &   0.6\% & 13.9\% & 0.18 &    6.0\% &  47.0\% & 0.09  \\
40 & 0.3 &    45.66 &   4.4\% & 84.3\% & 0.35 &    8.3\% &  84.3\% & 0.20  \\ \hline
60 & 0.1 &    30.62 &   2.3\% & 31.2\% & 0.35 &   11.1\% &  63.7\% & 0.18  \\
60 & 0.2 &   160.34 &   0.0\% &  0.1\% & 1.15 &    3.8\% &  56.2\% & 0.43  \\
60 & 0.3 &   328.37 &   0.0\% &  0.0\% & 2.14 &    0.7\% &  15.3\% & 0.83  \\ \hline
80 & 0.1 &   244.91 &   0.0\% &  0.1\% & 1.90 &    2.0\% &  24.1\% & 0.49  \\
80 & 0.2 &   886.53 &   3.8\% & 32.5\% & 4.83 &    9.7\% &  32.5\% & 1.33  \\
80 & 0.3 &  2022.59 &   0.2\% &  4.2\% & 8.15 &    1.7\% &  25.1\% & 2.43  \\ \hline
100 & 0.1 & 1694.93 &   0.5\% & 18.7\% &  4.12 &   3.7\% &  43.4\% & 1.00  \\
100 & 0.2 & 3203.70 &   0.5\% &  6.1\% & 11.51 &   1.3\% &  25.1\% & 2.91  \\
100 & 0.3 & 7122.23 &   0.4\% &  8.2\% & 21.96 &   0.7\% &   8.2\% & 5.58
\end{tabular}

\label{tab:compexp:rand:density}
\end{table}

As expected, the run time increases with an increasing number of nodes and edge density (the latter is due to an increase in the number of partition and identical coverage points, see Table~\ref{tab:compexp:rand:ppicp} for more details). While the average percentage deviation for the node-restricted algorithm is quite small, the maximal percentage deviation for an instance can be very large, indicating significantly inferior solutions; and even more so for the deterministic algorithm. In general, the average percentage deviation for the deterministic problems seems to be decreasing with an increase in the edge density $p$. This becomes more apparent once we average for each edge density $p$ over the number of nodes, obtaining average deviations of 9.63\%, 5.19\%, and 2.84\% for $p=0.1,0.2$, and $0.3$, respectively. Concerning the node-restricted problem, the opposite seems to be the case. Averaging the average percentage deviations over the number of nodes, we obtain percentages of 0.89\%, 1.23\%, and 1.23\%. Finally, if we average the average percentage deviations for each number of nodes over the edge densities, we can observe that the averages generally decrease with an increasing number of nodes. For the node-restricted algorithm (deterministic algorithm), we obtain averages of 1.9\%, 0.8\%, 1.3\%, and 0.4\% (12.0\%, 5.2\%, 4.5\%, and 1.9\%).

In Table~\ref{tab:compexp:rand:ppicp} we present the average number of partition points, $\#PP$, and identical coverage points, $\#ICP$ (which are not at the same time also partition points), as well as the average times to calculate them.

\begin{table}[htb]

\caption{Statistics for the number for partition and identical coverage points.}
\centering
\begin{tabular}{cc|cc|cc}
Nodes & Density &  \#PP & Time & \#ICP & Time \\ \hline
40 & 0.1 &     983.7 & 0.0 & 1018.1 &    0.6 \\
40 & 0.2 &    2119.7 & 0.0 & 1274.4 &    5.3 \\
40 & 0.3 &    2773.5 & 0.0 & 1670.9 &   16.3 \\ \hline
60 & 0.1 &    2735.1 & 0.0 & 1672.5 &    9.8 \\
60 & 0.2 &    4603.5 & 0.0 & 2375.3 &   53.6 \\
60 & 0.3 &    5824.7 & 0.0 & 2374.3 &  100.9 \\ \hline
80 & 0.1 &    4971.6 & 0.0 & 2712.2 &   81.6 \\
80 & 0.2 &    8066.8 & 0.0 & 3567.7 &  347.4 \\
80 & 0.3 &   10089.6 & 0.0 & 4228.7 &  713.8 \\ \hline
100 & 0.1 &   7622.7 & 0.0 & 3924.7 &  237.6 \\
100 & 0.2 &  11770.1 & 0.0 & 4979.5 &  906.8 \\
100 & 0.3 & 14640.46 & 0.0 & 6222.7 & 2423.6
\end{tabular}
\label{tab:compexp:rand:ppicp}
\end{table}

Next, we turn to analysing the effect of varying coverage radii. The results are given in Table~\ref{tab:compexp:rand:radii}, where $\%R$ denotes the percentage of the diameter of the graph that constitutes the coverage radius $R$, e.g., $\%R=0.1$ means that $R=0.1\cdot d^{max}$. All other columns are identical to Table~\ref{tab:compexp:rand:density}.

\begin{table}[htb]
\caption{Comparison for different coverage radii for the random graphs.}
\centering
\small
\begin{tabular}{cc|c|ccc|ccc}
& \multicolumn{1}{c}{} & \multicolumn{1}{c}{Max regret} & \multicolumn{3}{c}{Node-restricted} & \multicolumn{3}{c}{Deterministic} \\
Nodes & \%R &   Time  &  Avg Dev & Max Dev & Time  &  Avg Dev & Max Dev & Time \\ \hline
40 & 0.1 &   10.06   &  0.7\% & 13.9\% & 0.15 &  23.2\% &318.2\% & 0.08 \\
40 & 0.2 &   18.28   &  0.5\% &  7.6\% & 0.18 &   5.9\% & 73.5\% & 0.12 \\       
40 & 0.3 &   36.18   &  4.6\% & 84.3\% & 0.24 &   6.9\% & 84.3\% & 0.13 \\ \hline
60 & 0.1 &   90.29   &  0.0\% &  1.4\% & 0.99 &   4.3\% & 63.7\% & 0.41 \\       
60 & 0.2 &   157.42  &  2.2\% & 31.2\% & 1.17 &   5.9\% & 45.6\% & 0.49 \\       
60 & 0.3 &   271.63  &  0.0\% &  0.6\% & 1.49 &   5.4\% & 56.2\% & 0.54 \\ \hline
80 & 0.1 &   544.87  &  0.2\% &  4.2\% & 4.49 &   5.8\% & 28.4\% & 1.23 \\       
80 & 0.2 &   940.82  &  2.3\% & 32.5\% & 4.80 &   4.6\% & 32.5\% & 1.37 \\       
80 & 0.3 &   1668.34 &  1.5\% & 19.6\% & 5.60 &   3.0\% & 19.6\% & 1.66 \\ \hline
100 & 0.1 &  1669.05 &  0.2\% &  2.9\% & 11.01 &  3.6\% & 43.4\% & 2.79 \\
100 & 0.2 &  4053.70 &  0.4\% &  6.1\% & 11.44 &  0.8\% &  6.8\% & 2.99 \\
100 & 0.3 &  6298.12 &  0.8\% & 18.7\% & 15.15 &  1.3\% & 18.9\% & 3.72
\end{tabular}

\label{tab:compexp:rand:radii}
\end{table}

As expected, the run time increases with an increasing coverage radius. While, again, the average percentage deviations for the node-restricted algorithm are quite small (in contrast to the deterministic algorithm), the maximal percentage deviations can be very large. 
Concerning varying coverage radii, the individual averages don't show a consistent trend. However, if we average the average percentage deviations for each coverage radius over the number of nodes, we can again observe opposing trends for the deterministic and the node-restricted algorithm. For the former, the averages strictly decrease with an increasing coverage radius (9.22\%, 4.29\%, and 4.15\%), while for the latter they strictly increase (0.27\%, 1.36\%, and 1.72\%).

Finally, we turn to the effect of varying $UB$ for the generation of the upper and lower bounds on the edges. The results are given in Table~\ref{tab:compexp:rand:bounds}, where $UB$ denotes the value of the parameter used for generating the upper and lower bounds on the demand on an edge. All other columns are the same as before.
\begin{table}[htb]
\caption{Comparison for different values for $UB$ for the random graphs.}
\centering
\small
\begin{tabular}{cc|c|ccc|ccc}
& \multicolumn{1}{c}{} & \multicolumn{1}{c}{Max regret} & \multicolumn{3}{c}{Node-restricted} & \multicolumn{3}{c}{Deterministic} \\
Nodes & UB &   Time  &  Avg Dev & Max Dev & Time  &  Avg Dev & Max Dev & Time \\ \hline
40 &  10 &    21.5 &   2.0\% & 53.9\% &  0.2 &   9.1\% & 140.2\% & 0.1 \\
40 &  50 &    21.4 &   2.2\% & 84.3\% &  0.2 &  12.1\% & 251.1\% & 0.1 \\
40 & 100 &    21.4 &   1.6\% & 50.2\% &  0.2 &  14.8\% & 318.2\% & 0.1 \\ \hline
60 &  10 &   171.8 &   0.8\% & 26.4\% &  1.2 &   7.1\% &  63.7\% & 0.4 \\
60 &  50 &   171.5 &   0.6\% & 28.1\% &  1.1 &   4.0\% &  56.2\% & 0.4 \\
60 & 100 &   175.9 &   0.8\% & 31.2\% &  1.2 &   4.5\% &  46.4\% & 0.5 \\ \hline
80 &  10 &  1045.1 &   1.5\% & 29.1\% &  5.0 &   4.4\% &  29.1\% & 1.4 \\
80 &  50 &  1041.3 &   1.1\% & 25.8\% &  4.9 &   4.5\% &  28.0\% & 1.4 \\
80 & 100 &  1067.4 &   1.3\% & 32.5\% &  4.8 &   4.5\% &  32.5\% & 1.4 \\ \hline
100 &  10 & 3148.2 &   0.4\% &  8.2\% & 12.7 &   2.8\% &  43.4\% & 3.2 \\
100 &  50 & 3255.5 &   0.3\% &  5.7\% & 12.5 &   1.0\% &  19.9\% & 3.1 \\
100 & 100 & 3172.8 &   0.6\% & 18.7\% & 12.3 &   1.8\% &  30.2\% & 3.1
\end{tabular}

\label{tab:compexp:rand:bounds}
\end{table}
We can make very similar observations as for the previous two comparisons. A striking difference, however, is that the maximal percentage deviations for the node-restricted algorithm are very high for all combinations of $n$ and $UB$. Moreover, no consistent trend in the average percentage deviations can be observed with respect to increasing values of $UB$. This time, also averaging over the number of nodes for each value of $UB$ does not reveal anything (for the deterministic algorithm, we obtain averages of 5.85\%, 5.4\%, and 6.41\% for $UB=10,50$, and $100$, respectively, and for the node-restricted problem we have 1.19\%, 1.06\%, and 1.1\%). One might have expected that larger values of $UB$ would mean a larger ``range of uncertainty'' resulting in higher percentage deviations. But this is only evident for the deterministic algorithm for the 40-node instances. Observe that while the range of uncertainty is much larger in absolute values, it does not change in relative values. 

\subsubsection{Real-world graphs}

The two real-world data sets have 132 edges and 213 edges (with 106 nodes and 143 nodes, respectively). We carry out an analogous analysis as for the random data sets, starting with the coverage radii. The results are shown in Table~\ref{tab:compexp:street:radii}, where the values are just averaged over $UB$ now.
\begin{table}[htb]
\caption{Comparison for different coverage radii for the street graphs.}
\centering
\small
\begin{tabular}{cc|c|ccc|ccc}
& \multicolumn{1}{c}{} & \multicolumn{1}{c}{Max regret} & \multicolumn{3}{c}{Node-restricted} & \multicolumn{3}{c}{Deterministic} \\
Edges & \%R &   Time  &  Avg Dev & Max Dev & Time  &  Avg Dev & Max Dev & Time \\ \hline
132 & 0.1 &   26.0 & 0.1\% & 1.2\% & 0.2 & 34.1\% & 72.6\% & 0.1  \\
132 & 0.2 &  106.6 & 1.1\% & 3.4\% & 0.4 &  8.9\% & 32.7\% & 0.1  \\
132 & 0.3 &  196.8 & 1.1\% & 2.9\% & 0.5 & 25.1\% & 45.0\% & 0.1  \\ \hline
213 & 0.1 &  536.5 & 1.6\% & 7.4\% & 1.3 & 31.3\% & 61.2\% & 0.2  \\ 213 & 0.2 & 1652.2 & 3.6\% & 8.4\% & 2.1 & 12.7\% & 20.3\% & 0.4 \\ 213 & 0.3 & 2954.9 & 12.4\% & 25.3\% & 2.9 & 80.0\% & 108.8\% & 0.5 \end{tabular}

\label{tab:compexp:street:radii}
\end{table}

We can make similar observations concerning the run times. For the node-restricted algorithm, there seems to be a trend with respect to increasing values of $R$ resulting in increasing average deviations, but more tests would be required to verify this observation.

In Table~\ref{tab:compexp:street:ppicp} we present the number of partition points and identical coverage points (which are not at the same time also partition point), as well as the average times to calculate them. As $UB$ does not affect the number of points, the values are not averages but the actual numbers for both street graphs.

\begin{table}[htb]

\caption{Statistics for the number for partition and identical coverage points for the street graphs.}
\centering
\begin{tabular}{cc|cc|cc}
Edges & \%R &  \#PP & Time & \#ICP & Time \\ \hline
132 & 0.1 &   1451 & 0.0 &  15251 &  0.2 \\
132 & 0.2 &   1975 & 0.0 &  33460 &  1.8 \\
132 & 0.3 &   2481 & 0.0 &  48034 &  5.7 \\ \hline
213 & 0.1 &   3936 & 0.0 &  63138 &  8.8 \\
213 & 0.2 &   6035 & 0.0 & 110884 & 41.6 \\
213 & 0.3 &   7008 & 0.0 & 128882 & 70.0
\end{tabular}
\label{tab:compexp:street:ppicp}
\end{table}

Finally, we turn to the effect of varying $UB$ for the generation of the upper and lower bounds on the edges. The results are given in Table~\ref{tab:compexp:street:bounds}, where the values are just averaged over $\%R$.
\begin{table}[htb]
\caption{Comparison for different values for $UB$ for the street graphs.}
\centering
\small
\begin{tabular}{cc|c|ccc|ccc}
& \multicolumn{1}{c}{} & \multicolumn{1}{c}{Max regret} & \multicolumn{3}{c}{Node-restricted} & \multicolumn{3}{c}{Deterministic} \\
Edges & UB &   Time  &  Avg Dev & Max Dev & Time  &  Avg Dev & Max Dev & Time \\ \hline
132 &  10 &  101.3 & 1.0\% & 3.4\% & 0.3 & 30.3\% & 71.1\% & 0.1  \\
132 &  50 &  108.7 & 0.7\% & 3.3\% & 0.3 & 23.7\% & 72.6\% & 0.1  \\       
132 & 100 &  119.3 & 0.7\% & 2.8\% & 0.4 & 14.2\% & 37.3\% & 0.1  \\ \hline
213 &  10 & 1607.8 & 5.7\% & 25.3\% & 2.2 & 33.3\% & 94.1\% & 0.4  \\       
213 &  50 & 1941.1 & 5.0\% & 23.9\% & 2.1 & 43.4\% & 99.6\% & 0.3  \\       
213 & 100 & 1594.7 & 6.9\% & 23.9\% & 2.0 & 47.3\% & 108.8\% & 0.4
\end{tabular}

\label{tab:compexp:street:bounds}
\end{table}

As for the random data sets, no clear trend in the percentage deviations can be observed with respect to varying values of $UB$. The maximum deviations again underline that the node-restricted and the deterministic algorithm may produce significantly inferior solutions.

In conclusion, the obtained results show the advantages of using our model when we want to minimize the maximal regret. These advantages have significant consequences on the  service to be located as well as on the city's performance. For example, if a defibrillator has to be placed in the city with 213 edges (where UB=50) the obtained results means that in the worst case scenario for the demand (taking into account different coverage radii), the solution of our model leaves on average 5.0\% (43.4\%) less of the uncovered population by the node-restricted solution.
These results are even more remarkable if we consider the particular value of the radius that provided the maximum difference, where the solution of our model leaves 23.9\% (99.6\%) less of the uncovered population by node-restricted (deterministic) solution.

%

\section{The Unknown Linear Demand Case}
\label{sec:LinearDemand}%

In this section, we consider the case of an unknown linear demand realisation bounded by known linear lower and upper bound functions.
Let $lb_e(t) = a_e^{lb} + b_e^{lb}\cdot t$, $ub_e(t) = a_e^{ub} + b_e^{ub}\cdot t$, and $w_e(t) = a_e^{w} + b_e^{w} \cdot t$.
Unfortunately, Proposition~\ref{prop:Constant:FDS} and Observation~\ref{obs:Constant:DemandIndependentFDS} no longer hold for non-constant demand functions, although we can still easily compute the optimal solution of $\max_{x \in G} g(x,w)$ for a given $w$ in $\mathcal{O}(m^2 \log\,m)$ time (\citet{BerKalKra16}). What, however, still works is to  discretise the domain $w$ over which we optimise to find the worst-case demand realisation.

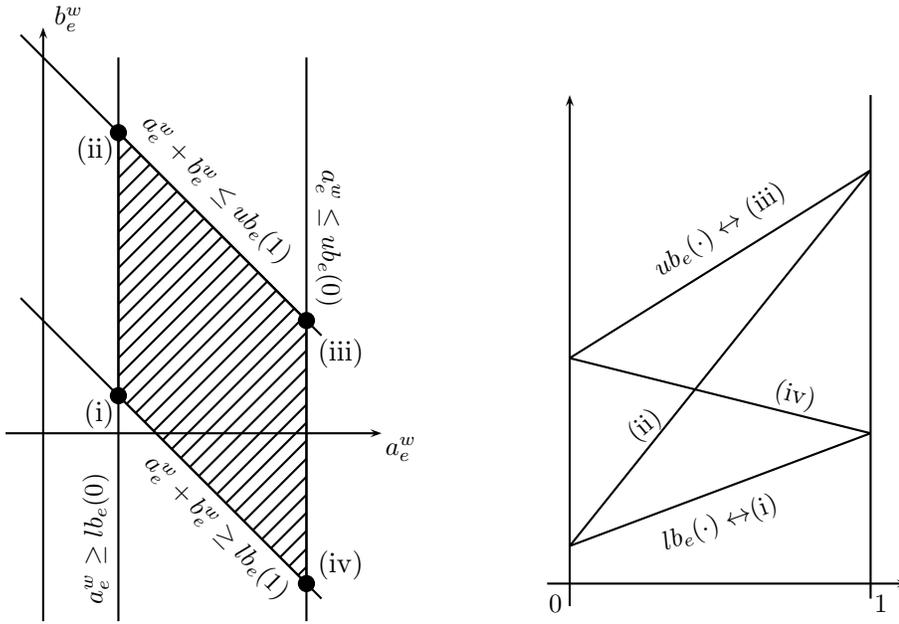
\begin{figure}[htb]
	\centering
	\begin{pspicture}(-1.0,-2.7)(12.7,5.4)
	
	\psline{->}(0.0,-2.5)(0.0,5.4) \psline{->}(-0.5,0)(4.5,0)
	\rput[tl](4.55,-0.05){\small $a_e^w$} \rput[bl](0.15,5.35){\small $b_e^w$}
	
	\psline(1.0,-2.5)(1.0,5.0) \rput{90}(0.7,-1.4){\small $a_e^w \ge lb_e(0)$}
	\psline(3.5,-2.5)(3.5,5.0) \rput{-90}(3.8,2.6){\small $a_e^w \le ub_e(0)$}
	\psline(-0.3,5.3)(3.7,1.3) \rput{-45}(2.3,3.25){\small $a_e^w + b_e^w \le ub_e(1)$}
	\psline(-0.3,1.8)(3.7,-2.2) \rput{-45}(2.3,-1.25){\small $a_e^w + b_e^w \ge lb_e(1)$}
	
	\pspolygon[fillstyle=hlines](1.0,0.5)(1.0,4.0)(3.5,1.5)(3.5,-2.0)
	\pscircle[fillstyle=solid,fillcolor=black](1.0,0.5){0.1} \rput[tr](0.95,0.45){(i)}
	\pscircle[fillstyle=solid,fillcolor=black](1.0,4.0){0.1} \rput[tr](0.95,3.95){(ii)}
	\pscircle[fillstyle=solid,fillcolor=black](3.5,1.5){0.1} \rput[tl](3.65,1.25){(iii)}
	\pscircle[fillstyle=solid,fillcolor=black](3.5,-2.0){0.1} \rput[bl](3.65,-1.95){(iv)}
	
	
	\psline{->}(7,-2.3)(7,4.5) \psline(11,-2.2)(11,4.5) 
	\psline{->}(6.7,-2.0)(11.5,-2.0) 
	\rput[tr](6.9,-2.15){\small $0$} \rput[tl](11.05,-2.15){\small $1$}
	
	\psline(7.0,-1.5)(11.0,0.0) \rput{20}(9,-1.2){\small $lb_e(\cdot) \leftrightarrow $(i)}
	\psline(7.0,1.0)(11.0,3.5) \rput{33}(9,2.7){\small $ub_e(\cdot) \leftrightarrow$ (iii)}
	\psline(7.0,-1.5)(11.0,3.5) \rput{40}(8,0.1){\small (ii)}
	\psline(7.0,1.0)(11.0,0.0) \rput{-15}(10,0.5){\small (iv)}
	\end{pspicture}
	\caption{Sketch of the feasible region $F_e$ (left-hand side) and the demand realisations corresponding to the four corners of $F_e$ (right-hand side).}
	\label{fig:Linear:FeasibleRegion}
\end{figure}

To that end, we first characterise the feasible region $lb \le w \le ub$ in terms of $a_e^{lb}, b_e^{lb}, a_e^{ub},$ and $b_e^{ub}$:
\begin{align}
\label{eq:Linear:FeasibleRegion}
 a_e^w & \:\ge\: lb_e(0) \,=\, a_e^{lb} \nonumber \\
 a_e^w & \:\le\: ub_e(0) \,=\, a_e^{ub} \nonumber \\
 lb \le w \le ub \qquad \Leftrightarrow \qquad a_e^w + b_e^w & \:\le\: ub_e(1) \,=\, a_e^{ub}+b_e^{ub} \\
 a_e^w + b_e^w & \:\ge\: lb_e(1) \,=\, a_e^{lb}+b_e^{lb} \nonumber \\
 a_e^w \,\ge\, 0, & \: b_e^w \in \mathbb{R} \nonumber
\end{align}
We denote by $F_e$ the feasible set of points $(a_e^{w},b_e^{w})$ satisfying the system of inequalities on the right hand side of \eqref{eq:Linear:FeasibleRegion}.
$F_e$ is a parallelogram in the $a_e^{w}/b_e^{w}-$space whose left and right  sides are vertical and whose upper and lower sides are diagonal with slope $-1$. See the sketch on the left-hand side in Figure~\ref{fig:Linear:FeasibleRegion}.

The four corners of the parallelogram are (i): $(a_e^{lb},b_e^{lb})$, (ii): $(a_e^{lb}, a_e^{ub}+ b_e^{ub}- a_e^{lb})$, (iii): $(a_e^{ub},b_e^{ub})$, and (iv): $(a_e^{ub},a_e^{lb}+b_e^{lb}-a_e^{ub})$.
The first and third coincide with $w_e(\cdot)=lb_e(\cdot)$ and $w_e(\cdot)=ub_e(\cdot),$ respectively. For the second and fourth, $w_e(\cdot)$ crosses diagonally between the lower and upper bound function. For example for the second point, $w_e(\cdot)$ starts at $lb_e(0)$ and ends at $ub_e(1)$.
See the sketch on the right-hand side in Figure~\ref{fig:Linear:FeasibleRegion}.
\medskip


\begin{prop}\label{prop:wcdemand}
	The worst-case demand realisation for a fixed $x,y \in G$ and $e\in E$ can be obtained by solving the following linear program: 
	\begin{eqnarray}\label{eq:Linear:MaxRegretEdgeFixedYNotInEdge}
	 r_e(x,y)=\max&& \hspace{-0.5cm} a_e^w \left(c_e(y) - c_e(x) \right)
	\,+\, \frac12 b_e^w\left({\bar{c}_e}(y) - {\bar{c}_e}(x) \right),   \\
	s.t.&&\hspace{-0.5cm} \eqref{eq:Linear:FeasibleRegion} \nonumber
	\end{eqnarray}   
 where 
	\[ {\bar{c}_e}(x) \:=\: \begin{cases} 1, & \text{ if } e \in E^c(x), \\
	1 - ((s_e^-(x))^2 - (s_e^+(x))^2), & \text{ if } e \in E^p(x), \\
	0, & \text{ if } e \in E^u(x), \\
	(s_{\ex}^-(x))^2 - (s_{\ex}^+(x))^2, & \text{ if } e = \ex.
	\end{cases}
	\]  
\end{prop} 

\ProofNoNL
Let $x = (\ex,\tx) \in G$, $y =(e_y,t_y) \in G$, and $e \in E$ be given. Next, we will derive the expression of the maximal regret for the different cases. We distinguish between the situations that define the different pieces of the functions $c_e$ and $\bar{c}_e$. We first consider the case that $e \in E^p(x) \cap E^p(y)$ and $x,y \not\in e$ (Case 1.1). Then, the maximal regret with respect to $x$ and $y$ on $e$ can be computed as
\begin{align} \label{eq:Linear:ProbMaxRegretEdgeFixedYNotInEdge}
r_{e}(x,y) & \:=\: \max_{lb_e \le w_e \le ub_e}\, (g_e(y,w_e)-g_e(x,w_e)) \nonumber \\
  & \:=\: \max_{lb_e \le w_e \le ub_e}\, \int_0^{s_e^+(y)} w_e(u)\, du + \int_{s_e^-(y)}^1 w_e(u)\, du \nonumber \\
  & \qquad \qquad \qquad \,-\, \int_0^{s_e^+(x)} w_e(u)\, du - \int_{s_e^-(x)}^1 w_e(u)\, du \nonumber \\
  & \:=\: \max_{lb_e \le w_e \le ub_e}\, \int_{s_e^+(x)}^{s_e^-(x)} w_e(u)\, du - \int_{s_e^+(y)}^{s_e^-(y)} w_e(u)\, du \nonumber \\
  & \:=\: \max_{lb_e \le w_e \le ub_e}\, \left[a_e^w\cdot u + \frac12b_e^w\cdot u^2 \right]_{s_e^+(x)}^{s_e^-(x)} - \left[a_e^w\cdot u + \frac12b_e^w\cdot u^2 \right]_{s_e^+(y)}^{s_e^-(y)} \nonumber \\
  & \:=\: \max_{lb_e \le w_e \le ub_e}\, a_e^w \left(s_e^-(x) - s_e^+(x) - s_e^-(y) + s_e^+(y) \right) \nonumber \\
  & \qquad \qquad \qquad \,+\, \frac12 b_e^w\left((s_e^-(x))^2 - (s_e^+(x))^2 - (s_e^-(y))^2 + (s_e^+(y))^2 \right) \nonumber \\
  & \:=\: \max_{lb_e \le w_e \le ub_e}\, a_e^w \left(c_e(y) - c_e(x) \right)
  \,+\, \frac12 b_e^w\left({\bar{c}_e}(y) - {\bar{c}_e}(x) \right). 
\end{align}

\medskip

So far we assumed that $x, y \not\in e$. Now, we assume that 
$y \not\in e$ but $x \in e$ (Case 1.2), then we obtain
\begin{align}
\label{eq:Linear:MaxRegretEdgeFixedYXInEdge}
r_e(x,y) 
  & \:=\: \max_{lb_e \le w_e \le ub_e}\, \int_0^{s_e^+(y)} w_e(u)\, du + \int_{s_e^-(y)}^1 w_e(u)\, du \,-\, \int_{s_e^+(x)}^{s_e^-(x)} w_e(u)\, du \nonumber \\
  & \:=\: \max_{lb_e \le w_e \le ub_e}\, \left[a_e^w\cdot u + \frac12b_e^w\cdot u^2 \right]_{0}^{s_e^+(y)} + \left[a_e^w\cdot u + \frac12b_e^w\cdot u^2 \right]_{s_e^-(y)}^{1} \nonumber \\
  & \qquad \qquad \qquad - \left[a_e^w\cdot u + \frac12b_e^w\cdot u^2 \right]_{s_e^+(x)}^{s_e^-(x)} \nonumber \\
  & \:=\: \max_{lb_e \le w_e \le ub_e}\, a_e^w \left(s_e^+(y) + 1 - s_e^-(y) - s_e^-(x) + s_e^+(x) \right) \nonumber \\
  & \qquad \qquad \qquad \,+\, \frac12 b_e^w\left((s_e^+(y))^2 + 1 - (s_e^-(y))^2 - (s_e^-(x))^2 + (s_e^+(x))^2 \right) \nonumber \\
  & \:=\: \max_{lb_e \le w_e \le ub_e}\, a_e^w \left(c_e(y) - c_e(x)\right)  \,+\, \frac12 b_e^w\left({\bar{c}_e}(y) - {\bar{c}_e}(x) \right).
\end{align}
Analogously, the same expression of $r_e(x,y)$ can be obtained in the remaining two subcases:  $x \not\in e, y \in e$ (Case 1.3) and $x, y \in e$ (Case 1.4). Next, we analyse 
the cases where $e \not\in E^p(x) \cap E^p(y)$. If $e \in E^p(y)\cap E^c(x)$ and $x,y \not\in e$ (Case 2.1), then
\begin{align}
\label{eq:Linear:MaxRegretEdgeFixedYNotInEdgeXFully}
r_e(x,y)
  & \:=\: \max_{lb_e \le w_e \le ub_e}\, \int_0^{s_e^+(y)} w_e(u)\, du + \int_{s_e^-(y)}^1 w_e(u)\, du \,-\, \int_{0}^{1} w_e(u)\, du \nonumber \\
  & \:=\: \max_{lb_e \le w_e \le ub_e}\, \left[a_e^w\cdot u + \frac12b_e^w\cdot u^2 \right]_{0}^{s_e^+(y)} + \left[a_e^w\cdot u + \frac12b_e^w\cdot u^2 \right]_{s_e^-(y)}^{1} \nonumber \\
  & \qquad \qquad \qquad - \left[a_e^w\cdot u + \frac12b_e^w\cdot u^2 \right]_{0}^{1} \nonumber \\
  & \:=\: \max_{lb_e \le w_e \le ub_e}\, a_e^w \left(s_e^+(y) + 1 - s_e^-(y) - 1 \right) \,+\, \frac12 b_e^w\left((s_e^+(y))^2 + 1 - (s_e^-(y))^2 - 1 \right) \nonumber \\
  & \:=\: \max_{lb_e \le w_e \le ub_e}\, a_e^w \left(c_e(y) - 1 \right)
  \,+\, \frac12 b_e^w\left({\bar{c}_e}(y) - 1 \right) \nonumber \\
  & \:=\: \max_{lb_e \le w_e \le ub_e}\, a_e^w \left(c_e(y) - c_e(x)\right)  \,+\, \frac12 b_e^w\left({\bar{c}_e}(y) - {\bar{c}_e}(x) \right).
\end{align}
 The remaining cases can be proven analogously. Therefore, we have proven that $r_e(x,y)$ can be obtained by solving a linear program.
\EndProofNoNL

	Since $r_e(x,y)$ can be reduced to solve a linear program with feasible region $F_e$, then, at least one of the four corners of $F_e$ is optimal.   In the following, we present the conditions  to identify which corner of $F_e$ will yield the maximal regret.

\begin{theorem}\label{prop:LinDemCas:OptDem}
		An optimal solution of \eqref{eq:Linear:MaxRegretEdgeFixedYNotInEdge}, $(a_e^{w*}, b_e^{w*})$, is given by the first column of Table~\ref{tab:OptDemRealisation} whenever the conditions of columns 2-4 are fulfilled. 
		
		\begin{table}[htbp]
		\caption{Worst-case demand realisation; optimal solution of \eqref{eq:Linear:MaxRegretEdgeFixedYNotInEdge}.}
			\centering
			\begin{tabular}{|l|c|c|c|}
				\hline
				\multicolumn{1}{|l|}{\multirow{2}{*}{$(a_e^{w*}, b_e^{w*})$}} & \multicolumn{3}{|c|}{Conditions}     \\ \cline{2-4}
				&   $c_e(y) - c_e(x)$    &   ${\bar{c}_e}(y) - {\bar{c}_e}(x)$    &      $({\bar{c}_e}(y) - {\bar{c}_e}(x)) -2(c_e(y) - c_e(x))$  \\
				\hline
				\multicolumn{1}{|l|}{$(a_e^{lb},b_e^{lb})$} &     $\leq0$  & $\leq0$      &   $\geq0$   \\
				\hline
				\multicolumn{1}{|l|}{\multirow{2}{*}{$(a_e^{lb}, a_e^{ub}+ b_e^{ub}- a_e^{lb})$}} & $\geq0$   &   $\geq0$     &    $\geq0$     \\ \cline{2-4}
				&   $\leq0$     &   $\geq0$     &      $-$   \\
				\hline
				\multicolumn{1}{|l|}{$(a_e^{ub},b_e^{ub})$} &   $\geq0$     &   $\geq0$     &      $\leq0$   \\
				\hline
				\multicolumn{1}{|l|}{\multirow{2}{*}{$(a_e^{ub},a_e^{lb}+b_e^{lb}-a_e^{ub})$}} &   $\leq0$     &     $\leq0$   &      $\leq0$    \\  \cline{2-4}
				&    $\geq0$    &    $\leq0$    &   $-$     \\
				\hline
			\end{tabular}%
			
			\label{tab:OptDemRealisation}%
		\end{table}%
		
	\end{theorem}
	
	\ProofNoNL
	The objective function in \eqref{eq:Linear:MaxRegretEdgeFixedYNotInEdge} is linear in $a_e^{w}$ and $b_e^{w}$. Hence, at least one of the four corners of $F_e$ will be optimal. Thus, if   $c_e(y) - c_e(x) \geq 0$ and $ {\bar{c}_e}(y) - {\bar{c}_e}(x) \leq 0$, an optimal solution is $(a_e^{ub},a_e^{lb}+b_e^{lb}-a_e^{ub}),$ i.e., corner (iv), since we want to make $a_e^w$ as large as possible and $b_e^w$ as small as possible to
		maximize the regret. Similarly, if  $c_e(y) - c_e(x) \leq 0$ and $ {\bar{c}_e}(y) - {\bar{c}_e}(x) \geq 0$, an optimal solution is $(a_e^{lb},a_e^{ub}+b_e^{ub}-a_e^{lb}),$ i.e., corner (ii), since we want to make $a_e^w$ as small as possible and $b_e^w$ as large as possible to
		maximize the regret. However,  if   $c_e(y) - c_e(x) > 0$ and $ {\bar{c}_e}(y) - {\bar{c}_e}(x) > 0,$ or    $c_e(y) - c_e(x) < 0$ and $ {\bar{c}_e}(y) - {\bar{c}_e}(x) < 0$, obtaining an optimal solution is not straightforward. Let us consider the case where $c_e(y) - c_e(x) > 0$ and $ {\bar{c}_e}(y) - {\bar{c}_e}(x) > 0,$ the other case can be analysed analogously. In this case, the two candidate points to be an optimal solution of \eqref{eq:Linear:MaxRegretEdgeFixedYNotInEdge} will be $(a_e^{lb}, a_e^{ub}+ b_e^{ub}- a_e^{lb})$ and $(a_e^{ub},  b_e^{ub}),$ i.e., corners (ii) and (iii). Thus, $(a_e^{ub},  b_e^{ub})$ will be an optimal solution of \eqref{eq:Linear:MaxRegretEdgeFixedYNotInEdge} if the following inequality holds:
	$$a_e^{lb}(c_e(y) - c_e(x))+\frac12(a_e^{ub}+b_e^{ub}-a_e^{lb})({\bar{c}_e}(y) - {\bar{c}_e}(x)) \leq a_e^{ub}(c_e(y) - c_e(x)) +\frac12 b_e^{ub}({\bar{c}_e}(y) - {\bar{c}_e}(x)),$$
	or  equivalently:
	$$\frac12(a_e^{ub}-a_e^{lb})({\bar{c}_e}(y) - {\bar{c}_e}(x)) \leq (a_e^{ub}-a_e^{lb})(c_e(y) - c_e(x)).$$
	By hypothesis $a_e^{ub}-a_e^{lb}\geq0.$ If $a_e^{ub}-a_e^{lb}=0$, corners (ii) and (iii) coincide, therefore we can assume without loss of generality that $a_e^{ub}-a_e^{lb}>0$. Thus, 
	$$\frac12({\bar{c}_e}(y) - {\bar{c}_e}(x)) \leq c_e(y) - c_e(x).$$ 
	It means that if $c_e(y) - c_e(x) > 0,$ $ {\bar{c}_e}(y) - {\bar{c}_e}(x) > 0,$ and ${\bar{c}_e}(y) - {\bar{c}_e}(x) \leq 2(c_e(y) - c_e(x))$, then an optimal solution of \eqref{eq:Linear:MaxRegretEdgeFixedYNotInEdge} is $(a_e^{ub},b_e^{ub})$, otherwise an optimal solution will be $(a_e^{lb}, a_e^{ub}+ b_e^{ub}- a_e^{lb})$. Therefore, we obtain the conditions described in Table \ref{tab:OptDemRealisation} and the result follows.
	\EndProofNoNL

\begin{theorem}\label{tm:LinDemCas:contquadratic}
For a given $e\in E, x\in \ex\in E,$ and $y\in \ey\in E,$ the conditions described in Theorem~\ref{prop:LinDemCas:OptDem} generate a subdivision of $\ex\times\ey$  with respect to $e$ with a constant number of cells, such that the representation of $r_e(x,y)$ over each cell is a quadratic function.
\end{theorem}
\ProofNoNL
Let $x \in \ex \in E$ and let $y \in \ey \in E$. From Lemma~\ref{lm:breakpoints} and Proposition~\ref{prop:ProbDescr:PartitionPoints} we know that the edge coverage functions, $s_{e}^+(x), s_{e}^-(x)$ $(s_{e}^+(y), s_{e}^-(y))$  are piecewise linear functions with a constant number of pieces over $\ex\in E$ $(\ey\in E)$. Thus, for each $e\in E$, the expressions $c_e(x)$ and $\bar{c}_e(x)$ $(c_e(y)$ and $\bar{c}_e(y))$ have a constant number of explicit  representations  for any $x \in \ex$ $(y \in \ey)$. 
 The breakpoints of the edge coverage functions for a given edge $e\in E$ correspond to either bottleneck points or network intersect points (Lemma \ref{lm:breakpoints}). Therefore, we add the corresponding vertical and horizontal lines induced by these constant number of partition points to the subdivision.  Moreover, by Theorem~\ref{prop:LinDemCas:OptDem} we can identify the corresponding worst-case demand realisations, i.e., $w_e^*$ in each cell. Therefore, a constant number of algebraic curves are derived from the conditions defined in Table~\ref{tab:OptDemRealisation}. The arrangement inside the square  $\ex\times \ey$ of these constant number of planar algebraic curves for a given $e\in E$, i.e., the vertical and horizontal lines induced by the breakpoints of the edge coverage functions and the conditions of Table~\ref{tab:OptDemRealisation},   generates a subdivision with a constant number of cells because each pair of curves intersects in a constant number of points. Furthermore, within each cell of this subdivision there is a common worst-case demand realisation, $w_e^*$, i.e., within each cell of this subdivision $r_e(x,y)$ has a unique representation as a quadratic function. 
\EndProofNoNL

In the following, we will compute the representation of the maximal regret in the network depicted in Figure 	\ref{fig:ExLD1} for a pair of edges. Moreover, we will represent the generated subdivision over the square composed by these two edges. 
\begin{example}\label{example} 
	Consider the network depicted in Figure \ref{fig:ExLD1}. For each edge $e\in E,$ its length is printed next to the edge. Let $R=1,$ $lb_{[1,2]}(t) = 3 -3t,$ $ub_{[1,2]}(t) = 15+7t,$ $lb_{[2,3]}(t) =3t,$ $ub_{[2,3]}(t) = 7+3t$,  $lb_{[1,3]}(t) = 2+3t$, and $ub_{[1,3]}(t) = 8+10t$. The set of partition points is given by $PP=V\cup \{([1,3],1/3), ([1,3],2/3),([2,3],1/2)\}.$ The three  partition points not included in $V$ indicated as dots in the figure are the  bottleneck point $([1,3],2/3)$ and the three network intersect points $(([1,3],1/3),([1,3],2/3),([2,3],1/2)$.
	
		The network of this example is the one used in Example \ref{example_cd}. Therefore, the edge coverage functions are identical. As stated before, the representation of these functions can be found in Appendix \ref{sec:app2}.

	Next, we denote $F_{[1,2]}, F_{[2,3]},$ and $F_{[1,3]}$ the feasible region $lb\leq w\leq ub$ in terms of $a_{e}^{lb},b_{e}^{lb}, a_{e}^{ub},$ and $b_{e}^{ub}$ for $e\in \{[1,2],[2,3],[1,3]\}$ respectively, i.e., the feasible set of points satisfying the system of inequalities of the right hand side of \eqref{eq:Linear:FeasibleRegion}.
	
	Let $x=([1,2], t_x) \in [1,2],$ i.e., $0\leq t_x\leq 1,$ and $y=([2,3],t_y)\in [2,3],$ i.e., $0\leq t_y\leq 1$. In the following, we compute the worst-case demand realisation in $(x,y)\in[1,2]\times[2,3]$ for a fixed edge  $[1,2] \in E$ applying Proposition \ref{prop:wcdemand}:
	$$r_{[1,2]}(x,y)= \begin{cases} \displaystyle \max_{ (a_{[1,2]}^w,\; b_{[1,2]}^w)\in F_{[1,2]}} \;-2 t_y a_{[1,2]}^w - 2 t_y^2 b_{[1,2]}^w, & \text{ if } 0\leq t_y \leq \frac{1}{2}, \\
	\displaystyle \max_{(a_{[1,2]}^w,\; b_{[1,2]}^w)\in F_{[1,2]}} \;-a_{[1,2]}^w-\frac12 b_{[1,2]}^w, & \text{ if } \frac{1}{2}\leq t_y \leq 1. \end{cases}$$ 
	
	First, we analyse the sign of the coefficients of $a_{[1,2]}^w$ and $b_{[1,2]}^w$ in both pieces. We observe that there is no change of sign in them, therefore the worst-case demand realisation is constant over each piece, i.e., there are only two cells generated by the horizontal line induced by the breakpoint of $r_{[1,2]}(x,y)$. From Proposition \ref{prop:LinDemCas:OptDem}, we obtain that the optimal solution is corner (i) for both pieces, thus:
	$$r_{[1,2]}(x,y)= \begin{cases}  -6 t_y + 6 t_y^2, & \text{ if } 0\leq t_y \leq \frac{1}{2}, \\
	-\frac{3}{2}, & \text{ if } \frac{1}{2}\leq t_y \leq 1. \end{cases}$$
	Figure \ref{fig:grafica0} shows the subdivision generated in the square $[1,2]\times [2,3],$ in each cell the expression of the maximal regret as a unique representation.
	\begin{figure}[h]
		\centering
		\includegraphics[height=0.35\linewidth]{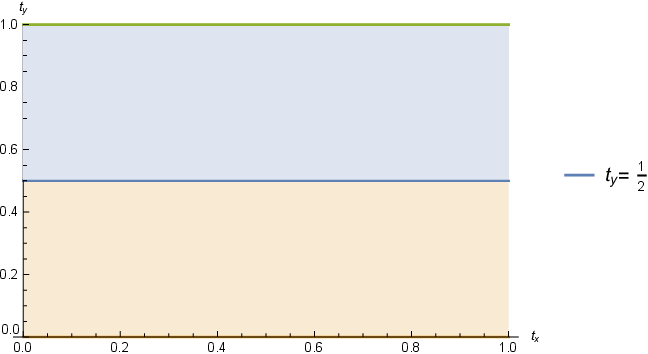}
		\caption{Cells for $r_{[1,2]}(x,y), \; (x,y)\in[1,2]\times[2,3]$.}
		\label{fig:grafica0}
	\end{figure}

Hereunder, we compute the worst-case demand realisation in $(x,y)\in[1,2]\times[2,3]$ for a fixed edge  $[2,3] \in E$ applying Proposition \ref{prop:wcdemand}:

	{\footnotesize
		$$r_{[2,3]}(x,y)= \begin{cases} \displaystyle \max_{(a_{[2,3]}^w,\; b_{[2,3]}^w) \in F_{[2,3]}} \;  \frac{1}{2}(1-t_x  + 2 t_y) a_{[2,3]}^w + \frac{1}{8} \left(-t_x^2 + (1 + 2 t_y)^2\right) b_{[2,3]}^w , & \text{ if } 0\leq t_y \leq \frac{1}{2}, \\
	\displaystyle \max_{ (a_{[2,3]}^w,\; b_{[2,3]}^w) \in F_{[2,3]}} \; \frac12 (3-2t_y-t_x) a_{[2,3]}^w +\frac12 \left(1-\frac{1}{4}\left((2t_y-1)^2+t_x^2\right)\right) b_{[2,3]}^w, & \text{ if } \frac{1}{2}\leq t_y \leq 1. \end{cases}$$}
	
	Since the coefficients of $a_{[2,3]}^w$  and $b_{[2,3]}^w$ have a non-negative value for both definitions, from Proposition \ref{prop:LinDemCas:OptDem}, we obtain that the changes in $r_{[2,3]}(x,y)$ are determined by the sign of $({\bar{c}_{[2,3]}}(y) - {\bar{c}_{[2,3]}}(x)) -2(c_{[2,3]}(y) - c_{[2,3]}(x))$. For $0\leq t_y \leq \frac{1}{2},$ this sign is also constant, so for this case new cells are not defined. However, for $\frac{1}{2}\leq t_y \leq 1,$ the sign is not constant, then we set $\left(1-\frac{1}{4}\left((2t_y-1)^2+t_x^2\right)\right) -(3-2t_y-t_x)=0$  and solve for $t_y\in\left[\frac12,1\right]$. We get the parametric curve $(t_x,\frac{3}{2} - \frac12 \sqrt{4t_x - t_x^2}),$ for  $2 - \sqrt{3} \leq t_x \leq 1$, this curve determines the change of the sign of  $({\bar{c}_{[2,3]}}(y) - {\bar{c}_{[2,3]}}(x)) -2(c_{[2,3]}(y) - c_{[2,3]}(x))$ and as a consequence, the change of definition of $r_{[2,3]}(x,y)$. 
	
	Figure \ref{fig:grafica1} shows the subdivision generated by the horizontal line induced by the breakpoint of $r_{[2,3]}(x,y)$ and the change of the sign of  $({\bar{c}_{[2,3]}}(y) - {\bar{c}_{[2,3]}}(x)) -2(c_{[2,3]}(y) - c_{[2,3]}(x))$. In each cell, the representation of $r_{[2,3]}(x,y)$ is unique: in the green area (upper left part) it is  $r_{[2,3]}(x,y)= \frac{7}{2}(3-2t_y-t_x)+\frac{3}{2}\left(1-\frac{1}{4}((2t_y-1)^2+t_x^2)\right)$; i.e., corner (iii) is the optimal solution, in the yellow one (upper right part) it is $r_{[2,3]}(x,y)= 5\left(1-\frac{1}{4}((2t_y-1)^2+t_x^2)\right);$ i.e., corner (ii) is the optimal solution, and in the blue one (lower part) it is 
	$r_{[2,3]}(x,y)= \frac{7}{2}(1-t_x  + 2 t_y) + \frac{3}{2} \left(-t_x^2 + (1 + 2 t_y)^2\right); $ i.e., corner (iii) is the optimal solution.
	
	\begin{figure}[h]
		\centering
		\includegraphics[height=0.35\linewidth]{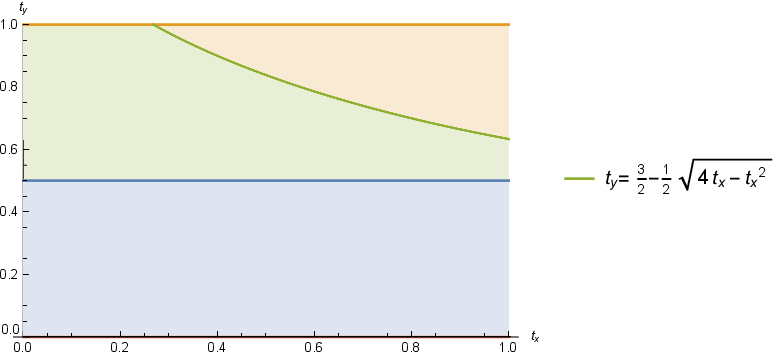}
		\caption{Cells for $r_{[2,3]}(x,y),$ $x=([1,2], t_x)$ and $y=([2,3],t_y).$ }
		\label{fig:grafica1}
	\end{figure}
	
	Next, we compute the worst-case demand realisation in $(x,y)\in[1,2]\times[2,3]$ for a fixed edge  $[1,3]\in E$ applying Proposition \ref{prop:wcdemand}:
	{\scriptsize $$r_{[1,3]}(x,y)= \begin{cases} \displaystyle \max_{ (a_{[1,3]}^w,\; b_{[1,3]}^w)\in F_{[1,3]}} \;  \frac{1}{3}(-1 + t_x)a_{[1,3]}^w -\frac{1}{18} (1 - t_x)^2b_{[1,3]}^w, & \text{if } 0\leq t_y \leq \frac{1}{2}, \\
		\displaystyle \max_{(a_{[1,3]}^w,\; b_{[1,3]}^w)\in F_{[1,3]}} \; \frac{1}{3}  (-2+t_x+2t_y) a_{[1,3]}^w+\frac12 \left(1-\frac{1}{9}\left((4-2t_y)^2+(1-t_x)^2\right)\right)b_{[1,3]}^w, & \text{if } \frac{1}{2}\leq t_y \leq 1. \end{cases}$$
	}
	
	From Proposition \ref{prop:LinDemCas:OptDem}, we obtain that the changes of definition of $r_{[1,3]}(x,y)$ are determined by the sign of  $c_{[1,3]}(y) - c_{[1,3]}(x)$, $\bar{c}_{[1,3]}(y) - \bar{c}_{[1,3]}(x),$ and  $(\bar{c}_{[1,3]}(y) - \bar{c}_{[1,3]}(x)) -2(c_{[1,3]}(y) - c_{[1,3]}(x)).$ In the Figure \ref{fig:grafica2}, the cells generated by the horizontal line induced by the breakpoint of $r_{[1,3]}(x,y)$ and the curves $c_{[1,3]}(y) - c_{[1,3]}(x)=0$, $\bar{c}_{[1,3]}(y) - \bar{c}_{[1,3]}(x)=0,$ and  $(\bar{c}_{[1,3]}(y) - \bar{c}_{[1,3]}(x)) -2(c_{[1,3]}(y) - c_{[1,3]}(x))=0$ are depicted. In each area, the representation of $r_{[1,3]}(x,y)$ is unique: in the pink area (lower part) it is $r_{[1,3]}(x,y)=\frac{2}{3}(-1 + t_x) -\frac{1}{6} (1 - t_x)^2;$ i.e., the corner (i) is the optimal solution, in the orange area (middle lower left part) it is $r_{[1,3]}(x,y)=\frac{2}{3}(t_x+2(t_y-1))+\frac{3}{2}\left(1-\frac{1}{9}((4-2t_y)^2+(1-t_x)^2)\right);$ i.e., the corner (i) is the optimal solution, in the blue and lavender area  (upper left and middle lower right part respectively) it  is $r_{[1,3]}(x,y)=\frac{2}{3}(t_x+2(t_y-1))+8\left(1-\frac{1}{9}((4-2t_y)^2+(1-t_x)^2)\right),$ i.e., the corner (ii) is the optimal solution, and in the green (upper right part) area it is $r_{[1,3]}(x,y)=\frac{8}{3}(t_x+2(t_y-1))+5\left(1-\frac{1}{9}((4-2t_y)^2+(1-t_x)^2)\right);$ i.e., the corner (iii) is the optimal solution. It is worth noting that although the expression $t_x=1-\frac{t_y}{2}$ implies a change in the sign of $c_{[1,3]}(y) - c_{[1,3]}(x)$, $r_{[1,3]}(x,y)$ does not change, because the optimal corner is the same. For this reason, the blue and the lavender cells would be considered as one in the following steps. 
	\begin{figure}[htb]
		\centering
		\includegraphics[width=0.7\linewidth]{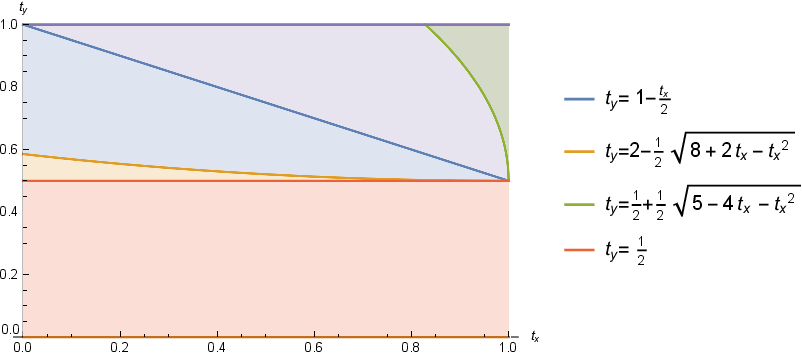}
		\caption{Cells for $r_{[1,3]}(x,y),$ $x=([1,2], t_x)$ and $y=([2,3],t_y).$ }
		\label{fig:grafica2}
	\end{figure}
	
	Combining all the above subdivisions, we obtain a finer subdivision shown in Figure~\ref{fig:grafica2a}. The function $r(x,y)$ has a unique  representation in each cell as a quadratic function  which is obtained as the sum of the corresponding $r_e(x,y),$ for each $e\in E.$
	
	\begin{figure}[htb]
		\centering
		\includegraphics[width=0.7\linewidth]{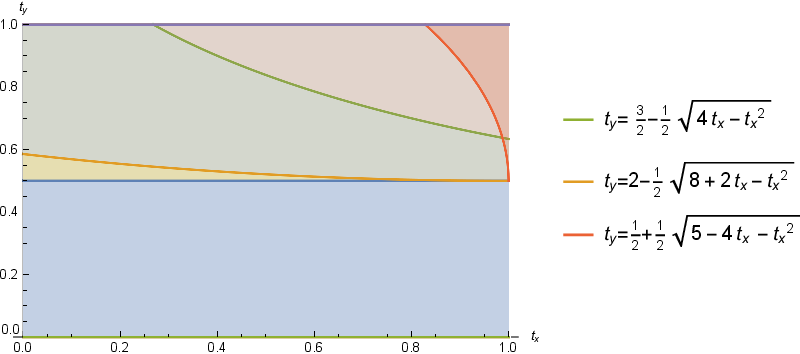}
		\caption{Cells for $r(x,y),$ $x=([1,2], t_x)$ and $y=([2,3],t_y).$ }
		\label{fig:grafica2a}
	\end{figure} 
\end{example}

\subsection{Resolution method}
In the current subsection, a polynomial time algorithm is developed for solving the single facility \name with unknown linear demand realisations. We start explaining the solution methods before analysing its complexity. The idea of the procedure is  to compute a subdivision over the square $\ex\times\ey$ for $\ex,\ey \in E$ where for any $x\in\ex,$ $\max_{y\in\ey}r(x,y)$ is achieved at the boundaries of the cells of this subdivision. We call $\mathcal{C}_{\ex \ey}$ the sets of arcs defining this subdivision.
\begin{theorem}
	The single facility \name with unknown linear demand realisations can be solved to optimality in polynomial time.
\end{theorem}
\ProofNoNL
Let $x \in \ex \in E$ and $y \in \ey \in E.$  For a given $e\in E$, we derive a subdivision over the square $\ex\times \ey$ with a constant number of cells, such that, within each cell  of the subdivision $r_e(x,y)$ has a unique representation (Theorem~\ref{tm:LinDemCas:contquadratic}). Indeed, within each cell of the arrangement generated by the vertical and horizontal lines induced by the breakpoints of the edge coverage functions and the algebraic curves defined on Table~\ref{tab:OptDemRealisation},  $r_e(x,y)$ has a unique representation  as a quadratic function for each $e\in E.$ Since $r(x,y)$ is the sum of $r_e(x,y)$ for all $e\in E$, the intersections of the cells generated for each $e \in E$ provide a finer subdivision  in the  square  $\ex\times \ey$ such that, within the new cells, $r(x,y)$ has a unique representation  as a quadratic function. Since $r(x,y)$ is quadratic,  for  a fixed $x$ the maximum of $r(x,y)$; i.e., $ \max_{y\in \ey} r(x,y)$ can be found on the boundary of the cells previously defined or, if the function is concave inside a cell, in the intersection of the curve $\frac{\partial r}{\partial y}(x,y)=0$ (for a fixed $x$) with the cell. Let $\mathcal{C}_{ \ex \ey}$ be  the set of the previous arcs (arcs defining the boundary of a cell and the arc $\frac{\partial r}{\partial y}(x,y)=0$ in the cells where $r(x,y)$ is concave) parametrised as $(x,y(x))$  in the square $\ex\times\ey$. The upper envelope, $h_{ \ex  \ey}(x)$, of $r(x,y(x))$ for all $(x,y(x))\in \mathcal{C}_{ \ex\ey}$  represents the optimal $y\in \ey$ for each $x\in \ex$ for the worst-case demand realisation (a similar idea of computing the upper envelopes over parametrised curves was used in \citet{LopPueRod13}); i.e.,
 \begin{equation}
h_{ \ex  \ey}(x) \:=\: \max_{lb \le w \le ub}\, \left( \max_{y\in \ey}\, g(y,w) \,-\, g(x,w) \right), \: \text{for } x \in \ex.
\end{equation}  
Repeating this procedure for each $ \ey\in E,$ the upper envelope, $h_{ \ex}(x)$, of 
 $r(x,y(x))$ for all $(x,y(x)) \in \bigcup\limits_{ \ey\in E}\mathcal{C}_{\ex \ey}$ determines the maximum regret of choosing $x \in \ex $ over the optimal location with respect to $w$ and $y\in G,$ i.e.,
 \begin{equation}
h_{ \ex}(x) \:=\: \max_{lb \le w \le ub}\, \left( \max_{y\in G}\, g(y,w) \,-\, g(x,w) \right), \: \text{for } x \in \ex.
\end{equation} 
Let  $x_{\ex}^{\prime}\in \ex$  be the  minimum  of $h_{\ex}(x)$ for $x \in \ex$. It determines the optimal location of the minmax regret problem in $\ex$; i.e., \begin{equation}
 h_{\ex}(x_{\ex}^{\prime}) \:=\: \min_{x \in \ex }\, \max_{lb \le w \le ub}\, \left( \max_{y\in G}\, g(y,w) \,-\, g(x,w) \right) .
\end{equation} 
This procedure should be repeated for each $\ex \in E$.
The optimal solution, $x^*,$ is $x^*=\displaystyle\min_{\ex \in E} h_{\ex}(x_{\ex}^{\prime}).$  Computing the upper envelope of these arcs and finding the minimum of this upper envelope can be done in polynomial time, therefore the result follows. The procedure presented in this proof is summed up in Algorithm \ref{algo:Linear}.
\EndProofNoNL
\begin{algorithm2e}[h]
	\DontPrintSemicolon \SetAlFnt{\small\sl}
	\SetAlCapFnt{\small\sl} \AlCapFnt
	
	\caption{Optimal algorithm for the single facility \name with unknown linear demand realisations.}
	\label{algo:Linear}
	\BlankLine 
	
	\KwIn{Network $G=(V,E)$; lower and upper bounds $lb_e(t) = a_e^{lb} + b_e^{lb}\cdot t$ and $ub_e(t) = a_e^{ub} + b_e^{ub}\cdot t$, respectively, for $e\in E$; coverage radius $R > 0$.}
	\BlankLine 
	\KwOut{Optimal solution $x^*$.}
	\BlankLine \BlankLine

	\nl \ForEach{$\ex\in E$}{
		
					\nl \ForEach{$\ey\in E$}{
			\nl Compute the subdivision generated by the arcs in $\mathcal{C}_{\ex \ey}:$ 
			\begin{enumerate}
				\item For each edge $e,$ the vertical and horizontal lines induced by the breakpoints of the edge coverage functions and the curves defining the conditions in Table \ref{tab:OptDemRealisation}.
				\item For any cell where $r(x,y)$ is concave in $y\in\ey$ for a fixed $x\in \ex$, the intersection of the curve $\frac{\partial r}{\partial y}(x,y)=0$ with that cell.
		\end{enumerate} 
			\BlankLine		
	}
\BlankLine
	\nl Obtain the upper envelope, $h_{\ex}(x),$ 
	of $r(x,y(x))$ of the arcs cointained in $\bigcup\limits_{\ey\in E}\mathcal{C}_{\ex \ey}$.  \;
			\nl Find the minimum $x_{\ex}^{\prime}$ of $h_{\ex}(x)$ over $\ex$.
			\BlankLine
			\nl \lIf{$h_{\ex}(x_{\ex}^{\prime}) < r(x^*)$}{set $x^* := x_{\ex}^{\prime}$, $r(x^*)= h_{\ex}(x_{\ex}^{\prime})$.} 
	}
	\BlankLine
	
	\nl \Return\ $x^*$.\;
\end{algorithm2e}

 For computing the complexity of solving the single facility \name with unknown linear demand realisations  on each edge by applying Algorithm~\ref{algo:Linear}, the following results are needed.

\begin{lemma}
	The set of arcs included in $\mathcal{C}_{\ex \ey}$ of the form $\frac{\partial r}{\partial y}(x,y)=0$ will be  constant functions.
\end{lemma}

\ProofNoNL
By definition, the function $r(x,y)$ does not have an $xy$-term. Therefore, $\frac{\partial r}{\partial y}(x,y)$ does not depend on $x$. Thus,  $\frac{\partial r}{\partial x \partial y}(x,y)=0$.
\EndProofNoNL

 In Algorithm~\ref{algo:Linear}, we are computing the upper envelope of $r(x,y(x)),$ where $(x,y(x))$ is the parametrisation of an algebraic arc of degree two, i.e., these curves may contain one square root and polynomials of degree two. Then, the intersection of two arcs could be transformed into a polynomial equation of degree eight, as a consequence they intersect at most eight times. Therefore, the complexity of the upper envelope of $m$ of these arcs is $\mathcal{O}(\lambda_{10}(m))$ and it can be computed in $\mathcal{O}(\lambda_{9}(m)\log m)$, see Theorem \ref{prop:compcostupperenvelope} in the Appendix.

Furthermore, the arrangement of $m$ planar algebraic curves inside a specified  square should be computed during the algorithm, this can be done in $\mathcal{O}(m^2)$ time, see \citet{KeyCulManKri00} for further details. Besides, this arrangement has complexity  $\mathcal{O}(m^2)$ as stated in \citet{HalMic17}.

\medskip
\noindent

\begin{theorem}
	The single facility \name with unknown linear demand realisations can be solved exactly in $\mathcal{O}(m\;\lambda_{10}(m^3))$ time using Algorithm~\ref{algo:Linear}.
\end{theorem}
\ProofNoNL  
Since $r(x,y)=\sum_{e \in E}r_e(x,y)$ and $r_e(x,y)$ has a constant number of representations  as quadratic functions over $x \in \ex$ and $y\in \ey,$ for each $e\in E$ (Theorem~\ref{tm:LinDemCas:contquadratic}), therefore $r(x,y)$ has $\mathcal{O}(m)$ representations  as quadratic functions over $x\in \ex, y\in \ey$. In the following, we will see the complexity of computing the expression of  $r(x,y)$. The expression of $r_e(x,y)$ is determined by the vertical and horizontal lines induced by the breakpoints of the edge coverage functions and the conditions of Table \ref{tab:OptDemRealisation}, i.e., for each $e\in E$ a constant number of algebraic curves subdivides the square  $\ex \times \ey$. Thus, the  square $\ex\times \ey$ is divided in $\mathcal{O}(m^2)$ cells by the $\mathcal{O}(m)$ algebraic curves previously obtained for all $e\in E$, this  arrangement can be constructed in $\mathcal{O}(m^2)$ time (\citet{KeyCulManKri00}). In a cell, we evaluate an interior point in order to apply Proposition \ref{prop:LinDemCas:OptDem} and identify $(a_e^w, b_e^w),$ for each $e\in E$, it is done in $\mathcal{O}(m)$ time. If we move to a neighbour cell (a cell that shares an arc boundary) only a constant number of addends changes, then we compute $r(x,y)$ in this  new cell in constant time, then we obtain $r(x,y)$ for all cells in  $\mathcal{O}(m^2)$ time just moving by adjacent cells in the square $\ex\times \ey$. 
Furthermore,  Step~3 includes in  $\mathcal{C}_{\ex \ey}$ $\mathcal{O}(m^2)$ arcs of cells boundaries and computes the derivative of  $\mathcal{O}(m^2)$ functions; it takes $\mathcal{O}(m^2)$ time. Step~4 obtain the upper envelope of  $\mathcal{O}(m^3)$ arcs; it takes  $\mathcal{O}(\lambda_{9}(m^3)\log m)$ (Theorem~\ref{prop:compcostupperenvelope}).  In addition, in Step~5, the minimum of an upper envelope is found, the complexity of it is the complexity of the upper envelope; i.e, $\mathcal{O}(\lambda_{10}(m^3))$ (Theorem~\ref{prop:compcostupperenvelope}).  
Finally, Step~6 takes constant time.

Therefore, since Step~4 and Step~5 are executed  once for each $\ex\in E$; i.e., $\mathcal{O}(m)$ times, the overall complexity of the algorithm is $\mathcal{O}(m\;\lambda_{10}(m^3))$.
\EndProofNoNL
\begin{observation}
	The complexity of Algorithm~\ref{algo:Linear} is  upper bounded by $\mathcal{O}(m^4 \log^{*} m)$,  see  Appendix.
\end{observation}

For illustrating the resolution method proposed in Algorithm \ref{algo:Linear}, we compute an iteration of the mentioned algorithm to solve the problem \name with unknown linear demand realisations in the network depicted in Figure \ref{fig:ExLD1}.

\ExampleCont{\ref{example}}
In Example \ref{example}, we computed the expression of the maximal regret of choosing $x\in[1,2]$ over $y\in[2,3],$ i.e., we computed the worst-case demand realisation. Now, the objective is to identify for each $x\in[1,2]$ the worst-case alternative for $y\in G$.

\begin{figure}[htb]
	\centering
	\includegraphics[width=0.7\linewidth]{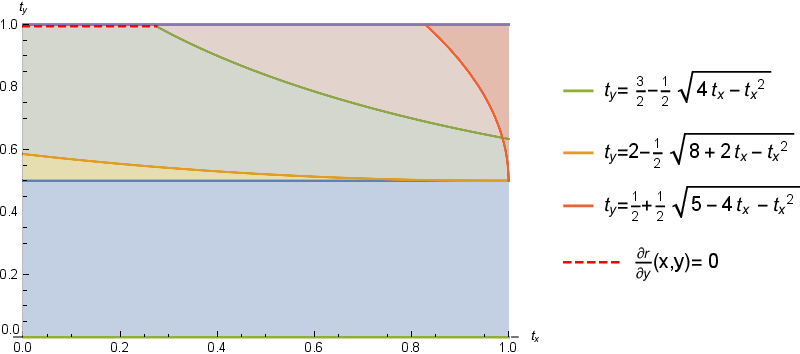}
	\caption{Cells for $r(x,y),$ $x=([1,2], t_x)$ and $y=([2,3],t_y).$ }
	\label{fig:grafica3v1}
\end{figure}

First, we consider the case $y\in[2,3]$. Since $r(x,y)$ is a quadratic function for a fixed $x\in[1,2]$, the $\max_{y\in[2,3]}r(x,y)$ will be found in the boundaries of the cells depicted in Figure \ref{fig:grafica2a}, or  in the curve $\frac{\partial r}{\partial y}(x,y)=0$  in the cells that $r(x,y)$ is concave. In Figure \ref{fig:grafica3v1}, these arcs are depicted. 

\begin{figure}[htb]
	\centering
	\includegraphics[width=0.7\linewidth]{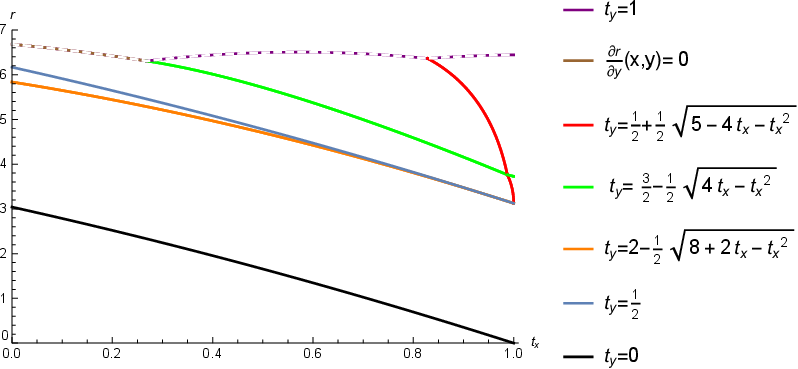}
	\caption{$r(x,y(x)),$ $x=([1,2], t_x)$ and $y=([2,3],t_y)$. }
	\label{fig:grafica4}
\end{figure}

The following step is to obtain $r(x,y(x))$ of the boundaries of the cells previously defined and the arcs $\frac{\partial r}{\partial y}(x,y)=0$ whenever $r(x,y(x))$  is concave. In Figure \ref{fig:grafica4} is depicted $r(x,y(x))$ of the previous mentioned arcs and the upper envelope is represented as a dotted function in the mentioned figure. Note that in the general algorithm there is no need to compute the upper envelope in this step, but it has been computed here in order to illustrate the algorithm and simplify the following graphic.

\begin{figure}[htb]
	\centering
	\includegraphics[width=0.7\linewidth]{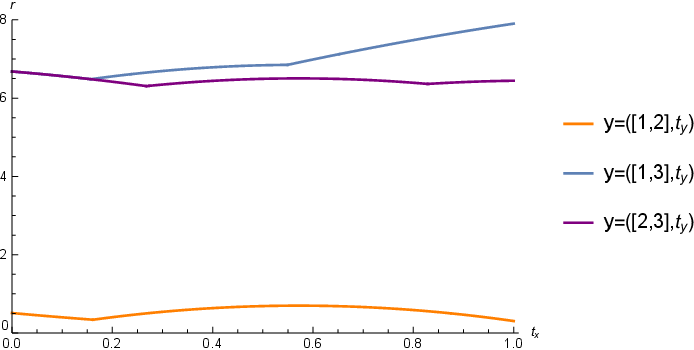}
	\caption{Upper envelope of $r(x,y(x)),$ $x=([1,2], t_x)$.}
	\label{fig:grafica5}
\end{figure}

Observe that for identifying the worst-case alternative the previous procedure should be repeated for each $y \in \ey\in E$, i.e., obtain the representation of $r(x,y)$, derive the subdivision of the square $\ex\times\ey,$ compute $r(x,y(x))$ of the boundaries of the cells, and calculate the upper envelope. The upper envelopes obtained for each $y \in \ey \in E$ are depicted in Figure \ref{fig:grafica5}. It worth mentioning that in the general algorithm, in this step we will compute the upper envelope of all the arcs previously obtained for each $\ey\in E$ instead of computing the upper envelope of the upper envelopes obtained for each edge, but in the example we did that in order to simplify the graphics. As can be appreciated, the minimum value of $r$ is $6.4836$, where $x_{[1,2]}^*=([1,2], 0.1572)$; these values were rounded to four decimal places.

This procedure should be repeated for each $x \in \ex\in E.$  The local minima over edge $[2,3]$ is vertex 2, where the value of $r$ is 7.9023, while the local minima over edge $[1,3]$ is the point $([1,3], 0.0533),$ where the value of $r$ is 6.3055; these values were rounded to four decimal places. The optimal location is the $x$ where the minimum value of $r(x^*)$ is found, i.e., the point $([1,3], 0.0533)$.

To highlight the usefulness of our max-regret approach, we also calculate the optimal location for a deterministic equivalent of the problem. We solve the problem assuming that the demand is known and it is equal to the mean of the upper and lower bound functions over each edge, i.e., $0.5(ub_e+lb_e).$ Following the procedure described in \citet{BerKalKra16}, the optimal solution is vertex 2, resulting in a maximal regret of $\dfrac{569}{72}\approx  7.9028$. Therefore, the difference between both optimal solutions concerning the max-regret value is significant. 
 The solutions of both models, the maximal regret of them and the covered demand assuming that the demand is deterministic are depicted in Table~\ref{tab:summary_2}.

\begin{table}[htbp]
 \caption{Solutions for Example~\ref{example}.}
  \centering
   {\tablinesep=2ex\tabcolsep=10pt \begin{tabular}{l|r|r|r}
    Model & Optimal solution & Maximal Regret & \textit{Covered Demand}  \\
    \hline
      \name     &  $\left([1,3],0.0533\right)$     &   $6.3055$    & 10.6858 \\
    \hline
       \citet{BerKalKra16}   &  $\left([1,2],1\right)$    &  $\dfrac{569}{72}$     & 12.1250 \\
    \end{tabular}%
    }
  \label{tab:summary_2}%
\end{table}%

\EndExample
\section{Conclusions and Outlook}
\label{sec:Conclusions}%
In this paper, we studied  the single-facility Minmax Regret Maximal Covering Location Problem \name on a network where the demand is unknown and distributed along the edges and the facility can be located anywhere on the network. We presented  two polynomial time algorithms for solving the cases where the realisation demands are unknown constant or linear functions. 

As far as we know, this is the first paper that applies the minmax
regret criterion on a maximal covering location problem on a network where the demand is distributed along the edges. Our results show that the problem is solvable in polynomial time (where the demand realisation are constant or linear functions) although the majority of polynomially solvable combinatorial optimization problems become NP-hard in the minmax regret version, as stated in \citet{kouYu97}.

There are several potential avenues for future research: Firstly, solving the single facility location problem for other kind of demand realisation functions. Secondly, considering the multi-facility location version of the problem. Finally, another interesting open question is how to  formulate the problem under a different criterion of coverage as e.g. the gradual covering, the cooperative covering model or the variable radius model, see \citet{BerDreKra10b} for more details of these criteria.

\vspace*{2ex}
\section*{Acknowledgements}
Thanks to the support of Agencia Estatal de Investigaci\'on (AEI) and the European Regional Development's funds (ERDF): project MTM2016-74983-C2-2-R, 2014-2020 ERDF Operational Programme and the Department of Economy, Knowledge, Business and University of the Regional Government of Andalusia: projects FEDER-UCA18-106895 and P18-FR-1422, Fundaci\'on BBVA: project NetmeetData (Ayudas Fundaci\'on BBVA a equipos de investigaci\'on cient\'ifica 2019), Telef\'onica and the BritishSpanish Society: Scholarship Programme 2018, and Universidad de C\'adiz: PhD grant UCA/REC01VI/2017. The authors would like to thank the anonymous reviewers for their comments and suggestions.

\bibliographystyle{abbrvnat}
\bibliography{./references}
\appendix
\section{Technical Notes}
In this appendix, we include some results about $\lambda_s(n)$, the maximum length of a Davenport-Schinzel sequence of order $s$ on $n$ symbols. We use these results for computing the complexity of the algorithms  proposed in the paper. 

\begin{theorem}[{{\citet[Theorem 6.5]{vsarir1995davenport}}}]\label{prop:compcostupperenvelope}
	Given a set of $m$ (unbounded $x-$monotone) Jordan arcs with at most s intersections between any pair  of arcs, its lower envelope has an $\mathcal{O}(\lambda_{s+2}(m))$ complexity, and it can be computed in $\mathcal{O}(\lambda_{s+1}(m) \log m)$ time.
\end{theorem}
 
Observe that $\lambda_1(m)=m,$ $\lambda_2(m)=2m-1,$ $\lambda_3(m)=\Theta(m\alpha (m)),$ and $\lambda_4(m)=\Theta(m2^{\alpha (m)}) $, where $\alpha (m)$ is the functional inverse of the Ackermann's function which grows very slowly. However, the problem of estimating $\lambda_s(m)$ for $s>4$ is more complicated.  For any constant $s$, it is well-known the bound $\lambda_s(m) = \mathcal{O}(m \log^* m)$. Recall that $\log^*m$ is the minimum number of times $q$ such that $q$ consecutive applications of
 the log operator will map $m$ to a value smaller than 1, i.e.,
 $\overbrace{ \log \dots \log m}  ^{(q)}\leq 1$. Actually, $\log^* m$ is the smallest height of
 an exponential ``tower'' of 2's, $2^{2^{2^{\dots}}}$
 which is $\geq m$ (nothing
 changes if $2$ is replaced by another base $b > 1$). Observe that,
 $\log^* m$ is much smaller than $\log m$ and it can be considered almost constant for ``practical'' values of $m$, see \citet{vsarir1995davenport} for further improvement of these bounds. 
 
\section{ Representation of the coverage functions of the examples}\label{sec:app2}
This appendix contains the representation of the edge coverage functions and the parts per unit of coverage functions of the network depicted in Figure \ref{fig:ExLD1}. These functions are used in Example \ref{example_cd} and \ref{example}.

Starting with edge $[1,2]$, for $x_1=([1,2], t_1)$ the edge coverage functions for all $e\in E$ are given by
	\begin{align*}
	& s_{[1,2]}^+(x_1)=0, & \quad & s_{[1,2]}^-(x_1
	)=1, \\
	&s_{[2,3]}^+(x_1)=\frac{t_1}{2},& &s_{[2,3]}^-(x_1)=1, \\
	&s_{[1,3]}^+(x_1)=\frac{1-t_1}{3},& &s_{[1,3]}^-(x_1)=1.
	\end{align*}
	Therefore, the parts per unit of coverage functions are: 
	\begin{align*}
	& c_{[1,2]}(x_1)=1, &  \\
	&c_{[2,3]}(x_1)=\frac{t_1}{2}, \\
	&c_{[1,3]}(x_1)=\frac{1-t_1}{3}.
	\end{align*}
	Likewise, for $x_2=([2,3], t_2)$ the edge coverage functions for all $e\in E$ are given by 
	\begin{align*}
	& s_{[1,2]}^+(x_2)=0, & \quad & s_{[1,2]}^-(x_2)=\begin{cases} 2t_2, & \text{ if } t_2 \leq \frac{1}{2}, \\
	1, & \text{ otherwise, } \end{cases}\\
	&s_{[2,3]}^+(x_2)=\begin{cases} \frac{2t_2-1}{2}, & \text{ if } t_2 \geq \frac{1}{2}, \\
	0, & \text{ otherwise } \end{cases} & & s_{[2,3]}^-(x_2)=\begin{cases} \frac{1+2t_2}{2}, & \text{ if } t_2 \leq \frac{1}{2}, \\
	1, & \text{ otherwise, } \end{cases} \\
	&s_{[1,3]}^+(x_2)=0, & & s_{[1,3]}^-(x_2)=\begin{cases} \frac{4-2t_2}{3}, & \text{ if } t_2 \geq \frac{1}{2}, \\
	1, & \text{ otherwise. } \end{cases}
	\end{align*}
	Thus, the parts per unit of coverage functions are: 
	\begin{align*}
	&  c_{[1,2]}(x_2)=\begin{cases} 1-2t_2, & \text{ if } 0 \leq t_2 \leq \frac{1}{2}, \\
	0, & \text{ if } \frac{1}{2} \leq t_2 \leq 1,  \end{cases}\\
	&  c_{[2,3]}(x_2)=\begin{cases} \frac{1}{2}+ t_2, & \text{ if } 0 \leq t_2 \leq \frac{1}{2}, \\
	\frac{3}{2}-t_2, & \text{ if } \frac{1}{2} \leq t_2 \leq 1,  \end{cases}\\
	&  c_{[1,3]}(x_2)=\begin{cases} 0, & \text{ if } 0 \leq t_2 \leq \frac{1}{2}, \\
	\frac{-1+2t_2}{3}, & \text{ if } \frac{1}{2} \leq t_2 \leq 1.  \end{cases}
	\end{align*}
	In addition, for $x_3=([1,3], t_3)$ the edge coverage functions for all $e\in E$ are given by 
	\begin{align*} 
	&s_{[1,2]}^+(x_3)=\begin{cases} 1-3t_3, & \text{ if } t_3 \leq \frac{1}{3}, \\
	0, & \text{ otherwise, } \end{cases}& &s_{[1,2]}^-(x_3)=1,\\
	&s_{[2,3]}^+(x_3)=0, & &s_{[2,3]}^-(x_3)=\begin{cases} \frac{4-3t_3}{2}, & \text{ if } t_3 \geq \frac{2}{3}, \\
	1, & \text{ otherwise, } \end{cases}\\
	&s_{[1,3]}^+(x_2)=\begin{cases} \frac{3t_3-1}{3} & \text{ if } t_3 \geq \frac{1}{3}, \\
	0 & \text{ otherwise, } \end{cases} & & 
	s_{[1,3]}^-(x_3)=\begin{cases} \frac{3t_3+1}{3}, & \text{ if } t_3 \leq \frac{2}{3}, \\
	1, & \text{ otherwise. } \end{cases}
	\end{align*}
	Hence, the parts per unit of coverage functions are: 
	\begin{align*} 
	&c_{[1,2]}(x_3)=\begin{cases} 1-3t_3, & \text{ if } 0\leq t_3 \leq \frac{1}{3}, \\
	0, & \text{ if } \frac{1}{3} \leq t_3  \leq 1, \end{cases}\\
	&c_{[2,3]}(x_3)=\begin{cases} 0, & \text{ if } 0\leq t_3  \leq \frac{2}{3}, \\
	-1+\frac{3t_3}{2}, & \text{ if } \frac{2}{3} \leq t_3 \leq 1, \end{cases}\\
	&c_{[1,3]}(x_3)=\begin{cases} \frac{1}{3}+t_3, & \text{ if } 0\leq t_3 \leq \frac{1}{3}, \\
	\frac{2}{3}, & \text{ if } \frac{1}{3} \leq t_3 \leq \frac{2}{3}, \\
	\frac{4}{3}-t_3, & \text{ if } \frac{2}{3} \leq t_3 \leq 1. \end{cases}
	\end{align*}

\end{document}